\documentclass[12pt]{article} %***
\usepackage[sectionbib]{natbib}
\usepackage{array,epsfig,fancyheadings,rotating}
\usepackage[]{hyperref}  %<----modified by Ivan
\usepackage{sectsty, secdot}
\sectionfont{\fontsize{12}{14pt plus.8pt minus .6pt}\selectfont}
\renewcommand{\theequation}{\thesection\arabic{equation}}
\subsectionfont{\fontsize{12}{14pt plus.8pt minus .6pt}\selectfont}
\usepackage{amsmath,amsfonts,bm}

\def\eqref#1{equation~\ref{#1}}
\def\1{\bm{1}}

\def\ry{{\textnormal{y}}}

\def\rvx{{\mathbf{x}}}
\def\rvy{{\mathbf{y}}}

\def\rmX{{\mathbf{X}}}

\def\vu{{\bm{u}}}
\def\vv{{\bm{v}}}

\def\vx{{\bm{x}}}

\def\vz{{\bm{z}}}

\def\mI{{\bm{I}}}

\DeclareMathAlphabet{\mathsfit}{\encodingdefault}{\sfdefault}{m}{sl}
\SetMathAlphabet{\mathsfit}{bold}{\encodingdefault}{\sfdefault}{bx}{n}

\newcommand{\R}{\mathbb{R}}

\def\R{\mathbb R}

\usepackage{graphicx}
\usepackage{subfigure}
\usepackage{amsmath,color}
\usepackage{algorithmicx,algorithm}
\usepackage{amssymb}
\usepackage{amsfonts}
\usepackage{multirow}
\usepackage{amsthm}
\numberwithin{equation}{section}
\newtheorem{theorem}{Theorem}%[section]
\newtheorem{Proposition}{Proposition}[section]

\newtheorem{Remark}{Remark}[section]
\usepackage{geometry}
\begin{document}
% 
% \tableofcontents\newpage

%%%%%%%%%%%%%%%%%%%%%%%%%%%%%%%%%%%%%%%%%%%%%%%%%%%%%%%%%%%%%%%%%%%%%%%%%%%%%%%%%%%%%%%%%%%%%%%%%%%%%%%%%%%%%%%%%%%%%%%%%%%%
%%%%%%%%%%%%%%%%%%%%%%%%%%%%%%%%%%%%%%%%%%%%%%%%%%%%%%%%%%%%%%%%%%%%%%%%%%%%%%%%%%%%%%%%%%%%%%%%%%%%%%%%%%%%%%%%%%%%%%%%%%%%

\renewcommand{\baselinestretch}{2}

% \markright{ \hbox{\footnotesize\rm Statistica Sinica
% %{\footnotesize\bf 24} (201?), 000-000
% }\hfill\\[-13pt]
% \hbox{\footnotesize\rm
% %\href{http://dx.doi.org/10.5705/ss.20??.???}{doi:http://dx.doi.org/10.5705/ss.20??.???}
% }\hfill }

% \markboth{\hfill{\footnotesize\rm FIRSTNAME1 LASTNAME1 AND FIRSTNAME2 LASTNAME2} \hfill}
% {\hfill {\footnotesize\rm FILL IN A SHORT RUNNING TITLE} \hfill}

% \renewcommand{\thefootnote}{}
% $\ $\par

%%%%%%%%%%%%%%%%%%%%%%%%%%%%%%%%%%%%%%%%%%%%%%%%%%%%%%%%%%%%%%%%%%%%%%%%%%%%%%%%%%%%%%%%%%%%%%%%%%%%%%%%%%%%%%%%%%%%%%%%%%%%

\fontsize{12}{14pt plus.8pt minus .6pt}\selectfont \vspace{0.8pc}
\centerline{\large\bf  Equivalence of state equations from different methods  }
\centerline{\large\bf  in High-dimensional Regression}
\vspace{2pt} 
%\centerline{\large\bf HERE IF A SECOND LINE IS NEEDED}
\vspace{.4cm} 
\renewcommand{\thefootnote}{\fnsymbol{footnote}}
\centerline{ Saidi Luo,  Songtao Tian\footnotemark[1]}
\footnotetext[1]{The authors contributed equally to this work and are listed in alphabetical order. Saidi Luo: First author, Center for Statistical Science, Tsinghua University. Songtao Tian: Corresponding author, Department of Mathematical Sciences, Tsinghua University.}
\vspace{.4cm} 
\centerline{\it Tsinghua University}
 \vspace{.55cm} \fontsize{9}{11.5pt plus.8pt minus.6pt}\selectfont

%%%%%%%%%%%%%%%%%%%%%%%%%%%%%%%%%%%%%%%%%%%%%%%%%%%%%%%%%%%%%%%%%%%%%%%%%%%%%%%%%%%%%%%%%%%%%%%%%%%%%%%%%%%%%%%%%%%%%%%%%%%%
\begin{quotation}
\noindent {\small Abstract:}\small ~State equations (SEs) were firstly introduced in the approximate message passing (AMP) to describe the mean square error (MSE) in compressed sensing. Since then a set of state equations have appeared in studies of logistic regression, robust estimator and other high-dimensional statistics problems. Recently, a convex Gaussian min-max theorem (CGMT) approach was proposed to study high-dimensional statistic problems accompanying with another set of different state equations. This paper provides a uniform viewpoint on these methods and shows the equivalence of their reduction forms, which causes that the resulting SEs are essentially equivalent and can be converted into the same expression through parameter  transformations. Combining these results, we show that these different state equations are derived from several equivalent reduction forms.
We believe that this equivalence will shed light on discovering a deeper structure in high-dimensional statistics.
% Functional sliced inverse regression (FSIR, \cite{ferre2003functional}) is the most recognized method in
% functional sufficient dimensional reduction (FSDR). However, there are still challenges in the application of FSIR. Among them, a remarking difficulty is choosing the number of slices in FSIR. It often depends on the sample size, the distribution of functional-valued variables, and other practical considerations. This greatly hinder the performance and development of FSIR in FSDR.  Aim to estimating the same central space as in FSIR, 
% we put forward a slicing-free inverse regression method in this paper. We achieve this by proposing the martingale difference divergence operator (MDDO) and corresponding algorithm. We establish a specific convergence rate of our method in estimating the central space.   Extensive simulations and real data experiments are shown to favor theoretic results and illustrate  the advantage of MDDO method.
% {\color{red}We extend the response from univariate to multivariate}

\vspace{9pt}
\par
\end{quotation}\par

	\def\thefigure{\arabic{figure}}
	\def\thetable{\arabic{table}}
	
	\renewcommand{\theequation}{\thesection.\arabic{equation}}

	\fontsize{12}{14pt plus.8pt minus .6pt}\selectfont
\newpage
\section{Introduction}
%\subsection{Motivations}
% In traditional statistical learning area, Maximum likelihood estimation(MLE) is often used to estimate model parameters for regression and classification problems.  In low-dimensional case where dimension $d$ is fixed and sample size  $n$ goes to infinite, a lot of precise and nearly perfect conclusions about MLE holds, such as asymptotic normality, asymptotically unbiased property and Wilks’ theorem can give the asymptotic distribution of the likelihood-ratio test(LRT). For regression problems, there are some extensions of MLE, such as M-estimators, which has similarly good properties, such as asymptotic normality, asymptotically unbiased property. \citep{huber1973robust} pointed out that the asymptotic distribution of M-estimator is Gaussian with 0 mean.

%A striking phenomenon in big-data era is that 
Classical statistical methods often failed in the high-dimensional data where the number of features is larger than the number of observed samples. Studies in high-dimensional data have attracted lots of attentions in past decades.  
A set of state equations (SEs) were first introduced in approximate message passing (AMP) algorithm in \citep{donoho2009message}  to precisely characterize the mean-square-error (MSE) and the phase transition phenomenon for true signal recovery in compressed sensing (CS). 
Since then, SEs, associated to certain AMP algorithm, have played indispensable role in various high-dimensional problems. For example, \citep{donoho2011noise} investigated the phase transition phenomenon and  the precise MSE of LASSO estimator; \citep{donoho2016high} studied the  variance of asymptotic distribution of M-estimator; \citep{huang2020asymptotic} provided a precise characterization of min-max MSE of $l_1$ penalized robust M-estimator and the corresponding phase transition phenomenon.

% \citep{donoho2011noise} utilize it to study the for LASSO regression; \citep{donoho2016high} for M-estimator;  \citep{huang2020asymptotic} for $l_1$ penalized robust M-estimator. More concrete and clear citation of their work}.

Though the SEs were first introduced through certain AMP type algorithms, researchers meet them in a variety of models through different methods. 
For example, the SEs appeared in \citep{el2013robust} when they performed the leaving-one-out (LOO) analysis of M-estimator in high dimensions. They showed that asymptotic normality, asymptotically unbiased property also hold as in the low dimension, nevertheless the variance of  asymptotic distribution of M-estimators  is higher.
\citep{sur2019modern} employed the similar idea to analyze the properties of MLE in logistic regression where the SEs were used to show that (1) asymptotically unbiased property does not hold; (2) variance of  asymptotic distribution increases; (3) likelihood ratio test is not distributed as chi-square.
%Strikingly, SEs were also appeared in another approach which is called Leaving-one-out (LOO) method to describe asymptotic behaviour for several models, such as  for Robust M-estimators and  for logistic regression.
SEs also appeared in another line of researches where Thrompoulidies et al. performed analysis of a family of high-dimensional problems through the Convex Gaussian min-max theorem (CGMT). More precisely,  \citep{thrampoulidis2018precise} characterized the MSE precisely for general regularized M-estimator problem in high-dimensions; \citep{salehi2019impact} established the correlation and MSE of the resulting estimator of regularized logistic regression;  \citep{deng2019model} showed the changing trend of 
MSE with the growth of features in support vector machine and logistic regression.
% {\color{red}Lastly, an insightful series of works \citep{barbier2019optimal,ricci2009cavity,moore2014cavity,krzakala2016statistical,coja2018information,mezard2003cavity,del2014cavity} have utilized the SEs (named as cavity method in statistical physics) as a ubiquitous tool when they studied the high-dimensional statistical problem through the perspective of statistical physics. Importantly, this tool has exhibited as a powerful weapon in applications of a lot of fields\citep{mezard2009information,obuchi2016cross,vuffray2014cavity, lesieur2015mmse,lesieur2016phase}.}

Lastly, an insightful series of works \citep{barbier2019optimal,ricci2009cavity,moore2014cavity,krzakala2016statistical,coja2018information,mezard2003cavity,del2014cavity} have utilized the SEs (named as cavity method in statistical physics) as a ubiquitous tool when they studied the high-dimensional statistical problem through the perspective of statistical physics. Importantly, this tool has exhibited as a powerful weapon in applications of a lot of fields\citep{mezard2009information,obuchi2016cross,vuffray2014cavity, lesieur2015mmse,lesieur2016phase}.
Though many papers have explicitly written down the corresponding state equations, none of them have shown that these sets of state equations are compatible. To the best of our knowledge, only \citep{deng2019model}  mentioned  there is another set of state equation  but without any comparison.

Although SEs were proved to be important in high-dimensional problems,  it is awkward that for one specific problem, the resulting SEs from AMP, CGMT and LOO are different. To be more clear, let us take a look at logistic regression.  
The SEs derived from CGMT (\ref{SE of LR from CGMT}) in \citep{deng2019model} are obviously different from  the SEs derived from LOO (\ref{SE of LR from LOO}) in \citep{sur2019modern}.
%Similar phenomenon appeared in all the aforementioned results. 
This is annoying, since the asymptotic performance for a specific high-dimensional problem should be unique no matter which method was used.
%To the best our knowledge, there is not any comparison of these methods and corresponding SEs is absent, which greatly hinders revealing the essence of these methods and discovering a deeper structure in high-dimensional statistics.

Therefore, we  are interested in the following questions: 

% {\color{red}
% \emph{Are SEs derived from different methods all equivalent in some sense? If so, from what viewpoint these methods are equivalent and are there more inner equivalence?}

% Among them, as the most direct, accessible, basic tool, equivalence of SEs is the basis of equivalence of methods and more inner equivalence.}

\emph{Are SEs derived from different methods all equivalent in some sense? If so, from what viewpoint these methods are equivalent and are there more inner equivalence?}

Among them, as the most direct, accessible, basic tool, equivalence of SEs is the basis of equivalence of methods and more inner equivalence.
%Inspired by widespread SEs and their significance in characterizing asymptotic performance of high-dimensional problems, we character the equivalence of AMP, CGMT, LOO heuristically and the equivalence of corresponding SEs precisely in this paper.

{\bf Our contributions. }
We successfully show that  for various high-dimensional problem, the different sets of SEs derived through different methods are actually equivalent to each other. More precisely, we construct the equivalence between different sets of SEs through explicit parameter transforms for  LASSO, M-estimator and logistic regression. These transformations are inspired by the statistical meanings of certain quantities appeared in the SEs. Moreover, we also provide a heuristic explanation on the relation between the different methods: AMP, CGMT and LOO.
% {\color{red}
% To the best of our knowledge, this is the first work to clearly clarify the equivalence among SEs derived from different methods and try to establish the equivalence of different methods.}
To the best of our knowledge, this is the first work to clearly clarify the equivalence among SEs derived from different methods and try to establish the equivalence of different methods.
% ,   is the basis of more insightful
% equivalence}. 

%derived from these methods are illustrated to be the same form by parameter transformations. 

%we prove the equivalence of these SEs with the assistance of certain invariants which possess statistical meanings. The paramet

%The key to the construction of the invariants is the statistical meaning of unknowns in SE.  We also illustrate several concrete examples to verify the equivalence of SE. Specifically, several SEs of significant high-dimensional optimization problems, such as 

%By investigating procedures of these methods, 

% Though the comparison of reduction form between CGMT and LOO is not discussed detailed in this paper, we believe their reduction forms are essential equivalent as AMP and CGMT.

% \subsection{Outlines}
{\bf Outlines.} 
In section \ref{Equivalence Examples}, we show that the SEs for M-estimator from AMP, LOO and CGMT are equivalent to each other. In section \ref{Equivalence CGMT AMP}, we show the equivalence of SEs derived from AMP and CGMT for another example and explain the essential reasons behind this equivalence. In Section \ref{Equivalence CGMT LOO}, we illustrate the similar work regarding the equivalence between CGMT and LOO.  Section \ref{Discussion and future directions} provides some discussions and future directions. Most proofs are deferred to the appendix.

%\subsection{Notations}
% \section{Preliminaries}
% \subsubsection{Notations}
{\bf Notations. }Let $\mathcal{N}(\boldsymbol{0}, \mI_d), \mathcal{N}(0, 1)$  denote the $d$-dimensional standard Gaussian  distribution and $1$-dimensional standard Gaussian distribution respectively. For a vector $\vx$, we denote $\|\vx\|_p$ as the $l_p$ norm of $\vx$. For an integer $n$ we denote $[n] $ as $\{1,\cdots,n\}$. 
% For a sequence of random variables $\{X_n\}$ that almost surely converge (converge in probability) to random variable $X$ we denote as $X_n \xrightarrow{a.s.}(\xrightarrow{P}) X$. 
We abbreviate independent and identically distributed to i.i.d.. For a function $f: \R \mapsto \R$, variable $x\in \R$ and $t > 0$, we denote the Moreau envelope associated with $f$ as
\begin{equation}
   M_f(x;t) := \min_{z\in \R} f(z) + \frac 1{2t}(x-z)^2,
\end{equation}
and the proximal operator, which is the solution of this minimization as \begin{equation}
   Prox_f(x;t) := \arg\min_{z\in \R} f(z) + \frac 1{2t}(x-z)^2.
\end{equation}
For multi-dimensional case $\vx = (x_1,\cdots,x_d)^T \in \R^d$, Moreau envelope and proximal operator are applied element-wisely: $M_f(\vx;t) := (M_f(x_i; t)) \in \R^d$ and $Prox_f(\vx;t) := (Prox_f(x_i;t)) \in \R^d$.

\section{An illustrative example}\label{Equivalence Examples}

%appeared in high-dimensional M-estimator illustrate high-dimensional M-estimator  as a representative example to show the equivalence of SEs from AMP CGMT and LOO. The equivalence is shown by doing parameter transformations to convert the corresponding SEs into the same form.

% For almost all common optimization problems, such as M-estimators, regularized M-estimators, logistic regression, the SEs derived from AMP and CGMT can be converted into the same form by parameter transformations. The same thing happens to LOO and CGMT.

Suppose that 
%Let $\{(\rvx_i,\ry_i)\}_{i = 1}^n$ be the data set where 
$\rvx_i\overset{\text{i.i.d. }}{\sim}   \mathcal{N}(0,\frac 1d \mI_d)$ and  $\ry_i \in \R$ satisfying that
\begin{equation}
 \ry_i = \rvx_i^T \boldsymbol{\beta}^* +\epsilon_i, \quad \text{for } i\in [n], \label{linearmodel}
\end{equation}
where $\epsilon_i$ are drawn i.i.d. from distribution $P_\epsilon$ with mean 0 and variance $\sigma_*^2$. We assume that the entries $\beta_i^*$ of  $\boldsymbol{\beta}^*$ are independently distributed as $\Pi$ which has finite second moment $r^2_* = \mathbb{E}_{\beta\sim\Pi} \beta^2$.

Let $\rho$ be a non-negative convex function. We are interested in the the Mean-squared-error (MSE) performance  $
  \lim_{n,p\rightarrow\infty}\frac1n\|\boldsymbol\beta-\boldsymbol{\beta}^*\|^2 
 $ of the  M-estimator:
\begin{equation}
\hat{\boldsymbol\beta}=\arg\min_{\boldsymbol\beta}\sum_{i=1}^n\rho(\ry_i-\rvx_i^T\boldsymbol\beta),
\label{M-estimator}
\end{equation}
when both $n$ and $d$ go to infinity satisfying that
$\lim_{n,d\to\infty} \frac dn = \kappa_* \in (0,\infty)$.

This problem first studied by \citep{el2013robust} where they showed that the MSE of $\hat{\beta}$ can be characterized by a set of SEs. More precisely, they proved the following proposition.

\vspace{2mm}
\begin{Proposition}\citep{el2013robust}
Given ratio $\kappa_*<1$. Consider the following system of nonlinear equations (SEs) regarding $(\tau_1,\gamma_1)$:
\begin{equation}
\label{SE of M-estimator from LOO}
\begin{aligned}
1-\kappa_*&=\mathbb E[\frac{\partial Prox_\rho}{\partial x}(W_1+\tau_1 Z_1;\lambda_1)],\\
\kappa_* \tau_1^2&:=\mathbb E[W_1+\tau_1 Z_1-Prox_\rho(W_1+\tau_1 Z_1;\lambda_1)]^2,
\end{aligned}    
\end{equation}
where $W_1\sim P_\epsilon,Z_1\sim \mathcal{N}(0,1)$ is independent of $W_1$. 
If this system of nonlinear equations possesses a unique solution $(\bar{\tau_{1}},\bar{\lambda_{1}})$, then the $\bar{\tau}_{1}$ is exactly the MSE of $\hat{\boldsymbol\beta}$ appeared in (\ref{M-estimator}).

% Given ratio $\kappa_*$, the SEs with unknowns $(\tau_2,\lambda_2)$ are:

\end{Proposition}

The $M$-estimator was also studied by \citep{donoho2016high} where they proved the following proposition.
\vspace{2mm}
\begin{Proposition}\citep{donoho2016high} 
Given ratio $\kappa_*<1$. Consider the following system of nonlinear equations (SEs) regarding $(\tau_2,\gamma_2)$:
\begin{equation}
\label{SE of M-estimator from AMP}
\begin{aligned}
\tau_2^2&=\frac1\kappa_*\lambda_2^2\mathbb E[\frac{\partial M_\rho}{\partial x}(W_2+\tau_2 Z_2;\lambda_2)]^2,\\
\kappa_*&=\lambda_2 \mathbb E[\frac{\partial^2 M_\rho}{\partial x^2}(W_2+\tau_2 Z_2;\lambda_2)];
\end{aligned}    
\end{equation}
where $W_2\sim P_\epsilon,Z_2\sim \mathcal{N}(0,1)$ is independent of $W_2$. 
If this system of nonlinear equations possesses a unique solution $(\bar{\tau_{2}},\bar{\lambda_{2}})$, then the $\bar{\tau}_{2}$ is exactly the MSE of $\hat{\boldsymbol\beta}$ appeared in (\ref{M-estimator}).

% Given ratio $\kappa_*<1$, then the MSE of (\ref{M-estimator}) is exactly $\bar{\tau}_2$ of the solution of the following system of nonlinear equations(SEs) with unknown $(\tau_2,\gamma_2)$  

% \begin{equation}
% \label{SE of M-estimator from AMP}
% \begin{aligned}

% \end{aligned}
% \end{equation}

% where $W_2\sim P_\epsilon$ and $Z_2\sim \mathcal{N}(0,1)$ is independent of $W_2$.whenever this SEs have an unique solution $(\bar\tau_2,\bar\gamma_2)$.
\end{Proposition}

Moreover, inspired by the work \citep{thrampoulidis2014tight}, we employ the CGMT techniques to study the $M$-estimator and show that the asymptotic MSE can be characterized by the the following SEs. To avoid unnecessary digression, we defer the detailed proof to the appendix.
\begin{Proposition}\label{prop:2.3}
Given ratio $\kappa_*<1$. Consider the following system of nonlinear equations (SEs) regarding $(\tau_3,\alpha, \mu)$:
\begin{equation}
\label{SE  from CGMT in M-estimator}
\begin{aligned}
0=&\frac{\alpha}{2}-\tau_3\sqrt{\kappa_*}-\frac{\alpha}{\mu^2}\mathbb E[\frac{\partial M_{\rho}}{\partial t}(W_3+\tau_3 Z_3;\alpha/\mu)],\\
0=&-\mu\sqrt{\kappa_*}+\mathbb E[Z_3\frac{\partial M_{\rho}}{\partial x}(W_3+\tau_3 Z_3;\alpha/\mu)],\\
0=&\frac{\mu}{2}+\frac{1}{\mu}\mathbb E[\frac{\partial M_{\rho}}{\partial t}(W_3+\tau_3 Z_3;\alpha/\mu)],
\end{aligned}    
\end{equation}
where $W_3\sim P_\epsilon, Z_3\sim \mathcal{N}(0,1)$ is independent of $W_3$. 
If this system of nonlinear equations possesses a unique solution $(\bar{\tau_{3}},\bar{\alpha},\bar\mu)$, then the $\bar{\tau}_{3}$ is exactly the MSE of $\hat{\boldsymbol\beta}$ appeared in (\ref{M-estimator}).
% Given ratio $\kappa_*\in(0,\infty)$, then the MSE of (\ref{M-estimator}) is exactly $\bar{\tau}_3$ of the solution of the following system of nonlinear equations(SEs) with unknown $(\tau_3,\mu,\alpha)$:
% % Given ratio $\kappa_*$, the SEs with unknowns $(\mu,\alpha,\tau_3)$ are:
% \begin{equation}
% \label{SE  from CGMT in M-estimator}
% \begin{aligned}
% 0=&\frac{\alpha}{2}-\tau_3\sqrt{\kappa_*}-\frac{\alpha}{\mu^2}\mathbb E[\frac{\partial M_{\rho}}{\partial t}(W_3+\tau_3 Z_3;\alpha/\mu)]\\
% 0=&-\mu\sqrt{\kappa_*}+\mathbb E[Z_3\frac{\partial M_{\rho}}{\partial x}(W_3+\tau_3 Z_3;\alpha/\mu)]\\
% 0=&\frac{\mu}{2}+\frac{1}{\mu}\mathbb E[\frac{\partial M_{\rho}}{\partial t}(W_3+\tau_3 Z_3;\alpha/\mu)]
% \end{aligned}.
% \end{equation}
% where $W_3\sim P_\epsilon, Z_3\sim \mathcal{N}(0,1)$ is independent of $W_3$.
% whenever this SEs have an unique solution $(\bar\tau_3,\bar\alpha,\bar\mu)$.
\end{Proposition}
% The proposition \ref{prop:2.3}
% We refer to the technique developed in \citep{thrampoulidis2014tight} when giving proof of this proposition.
\noindent The proof of this proposition is deferred to the Appendix S1.

On the one hand, these three sets of SEs are different at the first glance. On the other hand, since they are all supposed to describe the MSE of the $M$-estimators in high dimension, there shall be some relation between these three sets of equations. 
A striking fact is that we can actually show that all these three set of SEs are equivalent to each other.  More precisely, we have the following theorem. 

\vspace{2mm}
\begin{theorem}\label{prop para trans in M-estimator}
For M-estimator(\ref{M-estimator}), the SEs derived from AMP (\ref{SE of M-estimator from AMP}), LOO (\ref{SE of M-estimator from LOO}) and CGMT (\ref{SE from CGMT in M-estimator}) are equivalent. Specifically,
(\ref{SE of M-estimator from AMP}) can be converted into the same form as  (\ref{SE of M-estimator from LOO}) after  the following parameter transformations:
\begin{equation}
\label{para transformation for M-estimator in AMP and LOO}
\begin{aligned}
\tau_1=\tau_2, \quad 
\lambda_1=\lambda_2.\\
\end{aligned} 
\end{equation}
(\ref{SE of M-estimator from AMP}) can be converted into the same form as  (\ref{SE from CGMT in M-estimator}) after  the following parameter transformations:
\begin{equation}
\label{para transformation for M-estimator in AMP and CGMT}
\begin{aligned}
\tau_1=\tau_3, \quad
\lambda_1=\frac{\alpha}{\mu}.\\
\end{aligned} 
\end{equation}

\end{theorem}

The equivalence of these three sets of SEs seems straightforward, however, it suggests us that all the three procedures: AMP, CGMT and LOO might be deeply entangled in some sense. This will be investigated in this manuscript. 

The proof of this theorem is deferred to the Appendix S2.
%As will be detailed discussed in the section \ref{Equivalence of SEs from different methods}, the parameter transformations in this proposition depends on two invariants unchanged during the application of AMP and CGMT, which are main tool for constructing  parameter transformations in the following propositions.
%The detailed proof is deferred  until the appendix.
% \begin{Proposition}
% The SEs derived from LOO (\ref{SE of M-estimator from LOO}) are equivalent to the SEs derived from CGMT in M-estimator (\ref{SE from CGMT in M-estimator}). Specifically, 

% \begin{Remark}
% {\color{red}Since parameters in different methods express different statistical
% meanings in their own ways, the parameter transformations can help us
% related the statistical meanings in different methods. }
% % , which can help us
% % establish the equivalence of three methods, which we will try to explain in the following section. 
% \end{Remark}

% \end{Proposition}

\section{General results}\label{General results}

In this section, we show that the aforementioned equivalence between different sets of SEs holds not only for M-estimator, but also for LASSO and logistic regression in high dimensions.
\subsection{Equivalence between  the SEs derived from  CGMT and AMP}\label{Equivalence CGMT AMP}

%We first illustrate LASSO as another representative example to show the equivalence of SEs derived from AMP and CGMT.
Let us consider the following optimization:
% In this section, we only illustrate several representative examples to show the equivalence  of SEs for simplicity.
% In this section, LASSO is served as an example to indicate the equivalence between AMP and CGMT. 
\begin{equation}\label{LASSO}
\min_{\boldsymbol\beta} \frac12\|\rvy-\rmX\boldsymbol{\beta}\|_2^2+\lambda_*\|\boldsymbol{\beta}\|_1.
\end{equation}
where $\ry_{i}=\rvx_{i}^{\tau}\boldsymbol \beta^{*}+\epsilon_{i}$ with  $\rvx_i\overset{\text{i.i.d. }}{\sim}   \mathcal{N}(0,\frac 1n \mI_d)$ , $r_*:=\lim_{n,p\rightarrow\infty} \frac{\|\boldsymbol{\beta}^*\|}{\sqrt n}$ and
 $\lambda_* \ge 0$ is the regularized parameter, $\epsilon_i$ are drawn i.i.d. from distribution $P_\epsilon$ with mean 0 and variance $\sigma_*^2$.

We are interested in the the Mean-squared-error(MSE) performance  $
  \lim_{n,p\rightarrow\infty}\frac1n\|\boldsymbol\beta-\boldsymbol{\beta}^*\|^2 
 $ of the  LASSO.
% \citep{donoho2011noise} and \citep{mousavi2018consistent},\citep{bayati2011lasso},\citep{miolane2018distribution},\citep{javanmard2018debiasing} 
 \citep{donoho2011noise}, \citep{mousavi2018consistent}, \citep{bayati2011lasso, miolane2018distribution, javanmard2018debiasing}
have utilized the AMP to study the asymptotic performance of the LASSO estimator. For our purpose, we briefly recall the results in \citep{mousavi2018consistent} below.

\begin{Proposition}\citep{mousavi2018consistent}
Given noise scale $\sigma_*^2$  and ratio $\kappa_*$, consider the following system of nonlinear equations (SEs) regarding $(\tau_1,\gamma_1)$:
% Then for LASSO, the SEs derived from AMP which appeared in \citep{donoho2009message} are:

\begin{equation}
\label{SE of LASSO from AMP}
\begin{aligned}
\tau_1^2&=\sigma_*^2+\kappa_* \mathbb E[\eta(\beta_1+\tau_1 Z_1;\lambda_*+\gamma_1)-\beta_1]^2,\\
\gamma_1&=\kappa_*(\gamma_1+\lambda_*)\mathbb E[\eta'(\beta_1+\tau_1 Z_1;\lambda_*+\gamma_1)]
\end{aligned}
\end{equation}
where $Z_1\sim \mathcal N(0,1)$ is a standard normal variable, $\beta_1 \sim \Pi$ is independent of $Z_1$, $\eta(\cdot;\cdot)$ is the soft threshold function:
$$
\eta(x;t):=\text{sign}(x)(|x|-t)_+,
$$
where $x_+$ means $\max\{x,0\}$ and 
$$
\text{sign}(x):=\left\{ \begin{aligned} 1,\quad & \text{if } x>0, \\ 0,\quad &  \text{if } x=0; \\ -1,\quad & \text{if } x<0.
\end{aligned}
\right.
$$
If this system of nonlinear equations possesses a unique solution $(\bar{\tau_{1}},\bar{\lambda_{1}})$, then the $\bar{\tau}_{1}$ is exactly the MSE of $\hat{\boldsymbol\beta}$ appeared in (\ref{LASSO}).
\end{Proposition}
% $$
% \left\{
% \begin{aligned}
% -\frac{\alpha}{\sigma\tau_2}+\theta-1+\frac{\alpha+\lambda}{\lambda+1}&=0\\
% -\frac{1}{2\tau_2}+\frac{\gamma_2^2\alpha^2}{2\sigma^2\tau_2}-\frac{\tau_2\kappa_*}{2} \mathbb E[(Prox_{\tilde{f}}(\gamma_2Z+\theta \beta;\lambda_*))^2]+\frac{\sigma}{\lambda+1}&=0\\
% \gamma_2^2-r^2-\sigma_*^2+\frac{2[(\alpha+\lambda)r^2+\lambda\sigma_*^2]}{\lambda+1}-\frac{(\alpha+\lambda)^2r^2+\sigma^2+\lambda^2\sigma_*^2}{(\lambda+1)^2}&=0\\
% r^2\alpha-\sigma\tau_2\kappa_* \mathbb E[\beta Prox_{\tilde{f}}(\gamma_2Z+\theta \beta;\lambda_*)]&=0\\
% \sigma\tau_2\kappa_* \mathbb E[\frac{\partial Prox_{\tilde{f}}(\gamma_2Z+\theta \beta;\lambda_*)}{\partial Z}]&=\lambda\\
% \frac{\sigma}{2\tau_2^2}+\frac{r^2\alpha^2}{2\sigma\tau_2^2}-\frac{\sigma\kappa_*}{2} \mathbb E[(Prox_{\tilde{f}}(rZ+\theta \beta;\lambda_*))^2]&=0
% \end{aligned}
% \right.
% $$
% \ref{aa}

Inspired by the sequence of work \citep{thrampoulidis2014tight,thrampoulidis2015regularized,thrampoulidis2018precise,salehi2019impact}, we apply the CGMT to study the asymptotic performance of the LASSO estimator appeared in (\ref{LASSO}) and find that it is characterized by the following set of SEs. 
\begin{Proposition}
% Given ratio $\kappa_*<1$. Consider the following system of nonlinear equations (SEs) regarding $(\tau_1,\gamma_1)$  :
% \begin{equation}
% \label{SE of M-estimator from LOO}
% \begin{aligned}
% 1-\kappa_*&=\mathbb E[\frac{\partial Prox_\rho}{\partial x}(W_1+\tau_1 Z_1;\lambda_1)]\\
% \kappa_* \tau_1^2&:=\mathbb E[W_1+\tau_1 Z_1-Prox_\rho(W_1+\tau_1 Z_1;\lambda_1)]^2
% \end{aligned}    
% \end{equation}
% where $W_1\sim P_\epsilon,Z_1\sim \mathcal{N}(0,1)$ is independent of $W_1$. 
% If this system of nonlinear equations possesses a unique solution $(\bar{\tau_{1}},\bar{\lambda_{1}})$, then the $\bar{\tau}_{1}$ is exactly the MSE of $\hat{\beta}$ appeared in (\ref{M-estimator}).

Given noise scale $\sigma_*^2$, signal strength $r_*^2$ in model (\ref{linearmodel}) and ratio $\kappa_*$, Consider the following system of nonlinear equations (SEs) regarding $(\alpha,\sigma,\tau_2,\theta,\lambda,\gamma_2)$:
\begin{equation}
\begin{aligned}
0&=-\frac{\alpha}{\sigma\tau_2}+\theta-1+\frac{\alpha+\lambda}{\lambda+1},\\
0&=-\frac{1}{2\tau_2}+\frac{r_*^2\kappa_*\alpha^2}{2\sigma^2\tau_2}-\frac{\tau_2\kappa_*}{2} \mathbb E[(Prox_{\tilde{f}}(\gamma_2Z_2+\theta \beta_2;\lambda_*))^2]+\frac{\sigma}{\lambda+1},\\
0&=\gamma_2^2-r_*^2\kappa_*-\sigma_*^2+\frac{2[(\alpha+\lambda)r_*^2\kappa_*+\lambda\sigma_*^2]}{\lambda+1}-\frac{(\alpha+\lambda)^2r_*^2\kappa_*+\sigma^2+\lambda^2\sigma_*^2}{(\lambda+1)^2},\\
0&=r_*^2\kappa_*\alpha-\sigma\tau_2\kappa_* \mathbb E[\beta_2 Prox_{\tilde{f}}(\gamma_2Z_2+\theta \beta_2;\lambda_*)],\\
\lambda&=\sigma\tau_2\kappa_* \mathbb E[\frac{\partial Prox_{\tilde{f}}(\gamma_2Z_2+\theta \beta_2;\lambda_*)}{\partial x}],\\
0&=\frac{\sigma}{2\tau_2^2}+\frac{r_*^2\kappa_*\alpha^2}{2\sigma\tau_2^2}-\frac{\sigma\kappa_*}{2} \mathbb E[(Prox_{\tilde{f}}(\gamma_2 Z_2+\theta \beta_2;\lambda_*))^2].
\end{aligned}\label{SE of LASSO from CGMT}
\end{equation}
where $Z_2\sim \mathcal N(0,1)$ is a standard normal variable, $\beta_2 \sim \Pi$ is independent of $Z_2$, $\tilde f(x):=|x|$.

If this system of nonlinear equations possesses a unique solution $(\bar\alpha,\bar\sigma,\bar\tau_2,\bar\theta,\bar\lambda,\bar\gamma_2)$, then the $\frac{\bar{\lambda}_2}{\bar\theta}$ is exactly the MSE of $\hat{\boldsymbol\beta}$ appeared in (\ref{LASSO}).
\end{Proposition}

The detailed proof is deferred until the Appendix S3. The following proposition illustrate the equivalence between these two sets of SEs.

\vspace{3mm}
\begin{theorem}\label{Th, para transformation for LASSO}
The SEs of LASSO derived from AMP (\ref{SE of LASSO from AMP}) are equivalent to the SEs derived from CGMT (\ref{SE of LASSO from CGMT}). Specifically,   (\ref{SE of LASSO from CGMT}) can be converted into the same form as (\ref{SE of LASSO from AMP}) after  the following parameter transformations:
\begin{equation}
\label{para transformation for LASSO}
\begin{aligned}
\tau_1=\frac{\gamma_2}{\theta},\quad 
\gamma_1=\frac{\lambda_*}{\theta}-\lambda_*.
\end{aligned}
\end{equation}

\end{theorem}
\noindent The detailed proof is deferred until Appendix S4.

\vspace{2mm}
We provided a heuristic explanation on the equivalence of the SEs derived from AMP and CGMT.
For the sake of the self-contentment,  we briefly review the procedures of how to derive SEs from AMP and CGMT respectively.

\paragraph{Deriving SEs from AMP.}
The derivation of SE from AMP can be divided into two stages: 
%(1) construct iteration algorithm to solve the optimization problem  (2) derive the SEs. 

(1) Constructing an iterative algorithm. 
% AMP firstly complex initial problem to transform it into the problem of solving the Bayesian posterior distribution. During this step, objective function is transformed into a probability distribution, then based on the corresponding factor graph of this distribution, the message passing(MP) algorithm can be designed for solving Bayesian posterior distribution. Secondly, MP is reduced based on the aim to solving optimization problem, during which several approximation procedures are used such as large system limit, large $\beta$ limit and the approximation of iteration. 
\begin{itemize}
  \item [1)]
  AMP first transform initial optimization problem into pursuing a Bayesian posterior distribution where objective function is transformed into a probability distribution. 
  \item [2)]
Based on the corresponding factor graph of this distribution, it invokes the message passing(MP) algorithm to compute the Bayesian posterior distribution.
\item [3)] 
The MP is then further approximated by some  large system limit, large $\beta$ limit and the approximation of iteration.

\end{itemize}

(2) The asymptotic behavior of AMP is then characterized by  the state evolution equations/SEs.

\paragraph{Deriving SEs from CGMT.}
The derivation of SE from CGMT can be divided into four steps.
\begin{itemize}
  \item [1)]
The initial optimization problem is transformed into a min-max form, which is called  the primary optimization (PO) problem.
  \item [2)]
CGMT perform a dimensionality reduction on PO and obtain the auxiliary optimization (AO) problem
\item [3)] 
AO is further simplified to an optimization problem only depending on  several scalar variables, which is called scalar optimization (SO) problem.
\item [4)] 
SEs are derived by finding first-order optimality conditions of the asymptotic version of SO.
\end{itemize}

\begin{Remark}
We find that AO can be viewed as a relaxation of PO in the sense that  the feasible region of AO is larger than that of PO. Concrete examples, such as M-estimator, Logistic regression, Support vector machine and so on, are deferred to the Appendix S6.  We believe that this relaxation can help us understand the equivalence between the resulting SEs from AMP and CGMT respectively. 
\end{Remark}

% We provide a uniform viewpoint on AMP and CGMT in this part. Since procedures of deriving SEs from AMP can be divided into two parts: one is to solvie the original optimization problem and the other is to derive the asymptotic results of the solution. The Former  corresponds to the first step of CGMT, which means that the iteration of AMP is actually equivalent to the process of solving PO, as we will show in the following proposition. In other words, the iteration of AMP contains sufficient information of initial optimizations problem. The latter corresponds to the other steps of CGMT. 

We now present a uniform viewpoint on AMP and CGMT: 1) Constructing the AMP  corresponds to the first step of CGMT in LASSO, which suggests that the iteration of AMP is actually equivalent to the process of solving PO. 2) Deriving the SEs from AMP  corresponds to the last three steps of CGMT. Both of them aim to deriving SEs and characterizing asymptotic performance   by approximating  the initial optimization problem.
We proved the first statement in Proposition \ref{equivalence of AMP and CGMT}.

\vspace{2mm}
\begin{Proposition}\label{equivalence of AMP and CGMT}\citep{rangan2016fixed}
For LASSO, the fixed point of AMP is just the solution of first-order optimality conditions of PO in CGMT.
\end{Proposition}

\begin{proof}
For  CGMT, by introducing $\vu$ to constrain $ \vu = \rmX\boldsymbol\beta$ and Lagrange vector $\vv$, the corresponding PO can be written as:
$$
\min_{\boldsymbol\beta,\vu} \max_{\vv} \frac 1{2}\|\vu\|_2^2 - {\rvy}^T\vu+\frac12\|\rvy\|_2^2 +\lambda_*
\|\boldsymbol\beta\|_1 +  \vv^T(\rmX\boldsymbol\beta - \vu).
$$
Consider the first-order optimality conditions of PO:
\begin{equation}
\label{KKT}
\left\{
\begin{aligned}
0&=\lambda_* \text{sign}(\boldsymbol\beta) +\rmX^T\vv,\\
0&=\vu-\rvy-\vv ,\\
0&=\vu - \rmX\boldsymbol\beta.
\end{aligned}
\right.
\end{equation}

Comparing above formulas in (\ref{KKT}) leads to
$$
\lambda_* \text{sign}(\boldsymbol\beta) + \rmX^T(\rmX\boldsymbol\beta-\rvy) = 0.
$$
For AMP algorithm, the iteration of  LASSO is
$$
\begin{aligned}
{\boldsymbol\beta}^{t+1} &= \eta({\boldsymbol\beta}^t + \rmX^T\vz^t; \lambda_*+\gamma^t) ,\\
\vz^t &= \rvy - \rmX{\boldsymbol\beta}^t + \kappa_* \vz^{t-1}\langle \frac\partial{\partial x}\eta({\boldsymbol\beta}^{t-1}+\rmX^T\vz^{t-1}; \lambda_*+\gamma^{t-1})\rangle,\\
\gamma^{t} &= \kappa_*(\lambda_*+\gamma^{t-1})\langle \frac\partial{\partial x}\eta({\boldsymbol\beta}^{t-1}+\rmX^T\vz^{t-1}; \lambda_*+\gamma^{t-1})\rangle
,\end{aligned}
$$
% where for vector $\bold{x} = (x_i)\in \R^d$, $\langle \bold{x} \rangle := \sum_{i=1}^d x_i$ is the entry-sum of $\bold x$.
where $\frac{\partial}{\partial x}$ acts component-wisely.
For some vector $\vx$, $\langle \vx \rangle := \sum_{i=1}^d x_i$ denotes the entry-sum of $\vx$.

% for vector $\bold{x} = (x_i)\in \R^d$, $\langle \bold{x} \rangle := \sum_{i=1}^d x_i$ is the entry-sum of $\bold x$.

% Since the convergence point is a fixed point, $(\beta^\infty,z^\infty,\gamma^\infty)$ satisfy equations:
The fixed point $({\boldsymbol\beta}^\infty,\vz^\infty,\gamma^\infty)$ satisfy the following equations:
\begin{subequations}
\begin{align}
0&=(\lambda_*+\gamma)\text{sign}(\boldsymbol\beta) + \boldsymbol\beta - (\boldsymbol\beta+\rmX^T\vz) \label{fixed point equations 1},\\
\vz &= \rvy - \rmX\boldsymbol\beta + \kappa_* \vz\cdot c \label{fixed point equations 2},\\
\gamma &=  \kappa_*(\lambda_*+\gamma)c,\label{fixed point equations 3}
\end{align}
\end{subequations}
where $c = c(\boldsymbol\beta,\vz,\gamma) = \langle \frac\partial{\partial x}\eta(\boldsymbol\beta+\rmX^T\vz; \lambda_*+\gamma)\rangle$ and (\ref{fixed point equations 1}) is given by the following property about the soft thresholding function:
$$
t\cdot\text{sign}(z)+ z-x = 0,
$$
for $z = \eta(x;t)$ and some scalar $x$.
% $(z,x,t)$ satisfy equation $t\cdot\text{sign}(z)+ z-x = 0$ if $z = \eta(x;t)$. 

Simplifying (\ref{fixed point equations 2}) and (\ref{fixed point equations 3})       leads to:
\begin{equation}
\label{temp1}
\begin{aligned}
\vz &= \frac {\rvy-\rmX\boldsymbol\beta}{1-\kappa_*c} ,\\
\gamma &= \frac {\kappa_*c\lambda_*}{1-\kappa_* c}.
\end{aligned}
\end{equation}
Comparing (\ref{temp1}) with (\ref{fixed point equations 1}) gives
$$
\frac{\lambda_*}{1-\kappa_*c} \text{sign}(\boldsymbol\beta) - \frac 1{1-\kappa_*c}\rmX^T(\rvy-\rmX\boldsymbol\beta) = 0.
$$
which finishes the proof.
\end{proof}
% \begin{Remark}
% The solution of PO $\hat{\boldsymbol{\beta}}$ needs to be discussed componentwise. Denote $\hat{\boldsymbol{\beta}} = (\hat{\beta_i})$, the equation $\partial_\beta = 0$ in \ref{KKT} holds for each entry $i$ that $\hat\beta_i \neq 0$. For the entry $i$ that $\hat{\beta}_i = 0$, the optimality from PO gives $(X^TX\hat{\boldsymbol{\beta}}-y)_i\le \lambda_*$,  which is still equivalent to the fixed-point condition in AMP. Hence the equivalence holds for all entries of $\hat{\boldsymbol{\beta}}$.
% \end{Remark}

\begin{Remark}
It needs to be discussed component-wisely according to whether each entry of the optimal ${\boldsymbol{\beta}}$ is $0$ or not. The above proof holds for the entries that $\beta_i \neq 0$. For $i$ such that ${\beta}_i = 0$, the optimality from PO gives $-\lambda_* +(\rmX^T(\rmX\beta-\rvy))_i <0$ and $\lambda_* + (\rmX^T(\rmX\beta - \rvy))_i>0$. This is equivalent to $|(\rmX^T(\rmX{\boldsymbol{\beta}}-\rvy))_i|\le \lambda_*$, where $(\rmX^T(\rmX{\boldsymbol{\beta}}-\rvy))_i$ denote the $i$-th entry of $\rmX^T(\rmX{\boldsymbol{\beta}}-\rvy)$. This is still equivalent to the fixed-point condition in AMP. Hence the equivalence holds for all entries of ${\boldsymbol{\beta}}$.
% The solution of PO $\hat{\boldsymbol{\beta}}$ needs to be discussed componentwise. Denote $\hat{\boldsymbol{\beta}} = (\hat{\beta_i})$, the equation $\partial_\beta = 0$ in \ref{KKT} holds for each entry $i$ that $\hat\beta_i \neq 0$. For the entry $i$ that $\hat{\beta}_i = 0$, the optimality from PO gives $(X^TX\hat{\boldsymbol{\beta}}-y)_i\le \lambda_*$,  which is still equivalent to the fixed-point condition in AMP. Hence the equivalence holds for all entries of $\hat{\boldsymbol{\beta}}$.
\end{Remark}

% \begin{remark}
% %
% \end{remark}

\subsection{Equivalence between the SEs derived from  CGMT and LOO}\label{Equivalence CGMT LOO}
%In this part, we explain the equivalence between CGMT and LOO. We first illustrate logistic regression as another representative example to show the equivalence  of SEs derived from CGMT and LOO.

Suppose that % $\{(\rvx_i,\ry_i)\}_{i=1}^n$ where
$\rvx_i\overset{\text{i.i.d. }}{\sim}   \mathcal{N}(0,\frac 1d \mI_d)$ and $y_i\in \{-1,1\}$ drawn from logistic model:
\begin{equation}
 \mathbb{P}(\ry_i = 1|\rvx_i) = \rho'(\rvx_i^T\boldsymbol{\beta}^*), \quad \text{for } i\in [n], \label{logisticmodel}
\end{equation}
where $ \rho(t) = \log(1+e^t)$. Each entry of $\boldsymbol{\beta}$ is independently distributed as $\Pi$ which has finite second moment $r_*^2 = \mathbb{E}_{\beta\sim\Pi} \beta^2$.

We are interested in the following optimization problem:
%Logistic regression optimizes the following problem:
\begin{equation}
  \hat{\boldsymbol\beta}=\arg\min_{\boldsymbol{\beta}}\frac1n\sum_{i=1}^n\ell(\ry_i\rvx_i^T\boldsymbol\beta) ,\label{logistic regression}
\end{equation}
where $\ell(t):=\log(1+e^{-t})$.
When the $\hat{\boldsymbol\beta}$ exists, we are interested in the the Mean-squared-error(MSE) performance  
$\lim_{n,p\rightarrow\infty}\frac1n\|\boldsymbol\beta-\boldsymbol{\beta}^*\|^2 $ of the  Logistic regression. 

Logistic regression in high dimensions have been studied recently by \citep{candes2020phase,mousavi2018consistent}, \citep{deng2019model}. The asymptotic MSE of $\hat{\beta}$ was characterized by the following two propositions. 

\begin{Proposition}\citep{sur2019modern}
% Given ratio $\kappa_*<1$. Consider the following system of nonlinear equations (SEs) regarding $(\tau_1,\gamma_1)$  :
% \begin{equation}
% \label{SE of M-estimator from LOO}
% \begin{aligned}
% 1-\kappa_*&=\mathbb E[\frac{\partial Prox_\rho}{\partial x}(W_1+\tau_1 Z_1;\lambda_1)]\\
% \kappa_* \tau_1^2&:=\mathbb E[W_1+\tau_1 Z_1-Prox_\rho(W_1+\tau_1 Z_1;\lambda_1)]^2
% \end{aligned}    
% \end{equation}
% where $W_1\sim P_\epsilon,Z_1\sim \mathcal{N}(0,1)$ is independent of $W_1$. 
% If this system of nonlinear equations possesses a unique solution $(\bar{\tau_{1}},\bar{\lambda_{1}})$, then the $\bar{\tau}_{1}$ is exactly the MSE of $\hat{\beta}$ appeared in (\ref{M-estimator}).
Given signal strength $r_*^2$ in logistic model (\ref{logisticmodel}) and ratio $\kappa_*$, Consider the following system of nonlinear equations (SEs) regarding $(\lambda_1,\alpha_1,\sigma)$:
\begin{equation}
\label{SE of LR from LOO}
\begin{aligned}
\alpha_1^2&=\frac{1}{\kappa_*^2}\mathbb E[2\rho'(Q_1)\left(\lambda_1\rho'(Prox_\rho(Q_2;\lambda_1))\right)^2],\\
0&=\mathbb E[\rho'(Q_1)Q_1\lambda_1\rho'(Prox_\rho(Q_2;\lambda_1))],\\
1-\kappa_*&=\mathbb E[\frac{2\rho'(Q_1)}{1+\lambda_1\rho''(Prox_\rho(Q_2;\lambda_1))}],
\end{aligned}
\end{equation}
where
$$
% (Q_1,Q_2)\sim \mathcal N(\bold 0;\left(\begin{array}{c}
% r^2 &-\alpha_1 r^2\\
% -\alpha_1 r^2 &\alpha_1^2 r^2+\kappa_*\sigma^2
% \end{array}\right))
(Q_1,Q_2)\sim \mathcal N\left(\bold 0;
\begin{bmatrix} 

    r_*^2 & -\sigma r_*^2 \\

     -\sigma r_*^2 & \sigma^2 r_*^2+\alpha_1^2\kappa_*

      \end{bmatrix}\right),
$$
 and $\rho(t):=\log(1+e^t)$.

If this system of nonlinear equations possesses a unique solution $(\bar{\lambda_{1}},\bar{\alpha_{1}}, \bar{\sigma})$, then the MSE of $\hat{\boldsymbol\beta}$ appeared in (\ref{logistic regression}) is $[(\bar\sigma-1)\mathbb E_{\beta\sim\Pi}\beta]^2+{\bar\alpha}^2$.
\end{Proposition}

\begin{Remark}
In \citep{sur2019modern}, it is assumed that $X_{i,j}\sim \mathcal{N}(0,\frac{1}{n}I_d)$ and $r_*^2=\kappa_*\mathbb{E}_{\beta\sim\Pi} \beta^2$, which is slightly different from the setting in this paper. However, this difference only leads to a constant change related to $\kappa_*$ in the final parameter transformations (\ref{para transformations of LR from LOO and CGMT}) and does not affect the equivalence of these two set of SE. 
\end{Remark}

\begin{Proposition}\citep{deng2019model}
Given signal strength $r_*^2$ in logistic model (\ref{logisticmodel}) and ratio $\kappa_*$, Consider the following system of nonlinear equations (SEs) regarding $(\lambda_2,\alpha_2,\mu)$:
\begin{equation}
\label{SE of LR from CGMT}
\begin{aligned}
0&=\mathbb E[V\ell'(Prox_\ell(\alpha_2 Z+\mu V;\lambda_2))],\\
\alpha_2^2\kappa_*&=\lambda_2^2\mathbb E[(\ell'(Prox_\ell(\alpha_2 Z+\mu V;\lambda_2)))^2],\\
\kappa_*&=\lambda_2\mathbb E[\frac{\ell''(Prox_\ell(\alpha_2 Z+\mu V;\lambda_2))}{1+\lambda\ell''(Prox_\ell(\alpha_2 Z+\mu V;\lambda_2))}],
\end{aligned}   
\end{equation}
where $Z\sim \mathcal N(0,1)$, $V=Z_1Y_{r_*}$, in which $Z_1 \sim \mathcal{N}(0,1) $ is  independent of $Z$ and $\ Y_{r_*}\sim Ber(\rho'(r_*Z_1))$. $Ber(p)$ denotes the Bernoulli distribution with probability $p$ for the value $+1$ and probability $1-p$ for the value $-1$. If this system of nonlinear equations possesses a unique solution $(\bar{\lambda_{2}},\bar{\alpha_{2}}, \bar{\mu})$, then the MSE of $\hat{\boldsymbol\beta}$ appeared in (\ref{logistic regression}) is $[(\frac{\bar{\alpha}_2}{\sqrt{\kappa_*}}-1)\mathbb E_{\beta\sim\Pi}\beta]^2+(\frac{\bar\mu}{r_*})^2$.
\end{Proposition}

As before, we can show that (\ref{SE of LR from LOO}) and (\ref{SE of LR from CGMT}) are equivalent. 
\begin{theorem}\label{Th,para transformations of LR from LOO and CGMT}
For logistic regression (\ref{logistic regression}), the SEs derived from LOO (\ref{SE of LR from LOO}) and CGMT (\ref{SE of LR from CGMT}) are equivalent. Specifically, (\ref{SE of LR from LOO}) can be converted into the same form as (\ref{SE of LR from CGMT}) after  the following parameter transformations:
\begin{equation}
\label{para transformations of LR from LOO and CGMT}
\begin{aligned}
\alpha_1=\frac{\alpha_2}{\sqrt{\kappa_*}},\quad 
\sigma&=\frac\mu{r_*},\quad
\lambda_1=\lambda_2.
\end{aligned}
\end{equation}

\end{theorem}
\noindent The proof of this theorem is deferred to Appendix S5.

% \subsubsection{Procedures of deriving SEs from LOO}

For the sake of self-contentment, we briefly review the procedure on deriving SEs from LOO. 

% \begin{Remark}
% {\color{red}The similarity of
% parameter transformations in Theorem \ref{prop para trans in M-estimator}, Theorem \ref{Th, para transformation for LASSO} and Theorem \ref{Th,para transformations of LR from LOO and CGMT} strongly suggest the equivalence among different methods since different methods are using similar statistics to express asymptotic properties.}
% \end{Remark}

\paragraph{Deriving SEs from LOO} The derivation can be divided into 4 steps.
% For the original optimization, this approach first consider its first-order condition and derive an equation for $\hat\beta$. Analogously, two equations can also be obtained by calculating the first-order condition after leaving out one observation and one predictor respectively. Subtracting these two equations from original equation and performing Taylor expansion, two approximate equations are derived which characterize the relationship between the estimator after leaving one predictor and one observation. After this process, the correlation between the estimator $\hat{\boldsymbol{\beta}}$ and sample $(\rmX, \rvy)$ can be basically eliminated, which means that the relationship between the original sample and the parameter estimator can be characterized by two independent random vectors. Hence the previous high-dimensional(dimension $d$) approximate equations can be relatively easily transformed in 1-dimensional equations when calculating the asymptotic statistics such as MSE and bias as $n,d\to \infty$. Each statistic corresponds to an equation, which can be combined to derive the set of state equation.

\begin{itemize}
  \item [1)]
First, for the original optimization problem,  LOO considers first-order conditions of three cases: a) keeping all observations and predictors, corresponding solution is denoted as $\hat\beta$, b) leaving one predictor, corresponding solution is denoted as $\hat\beta_{(-j)}$ and c) leaving one predictor and one observation, corresponding solution is denoted as $\hat\beta_{(-i),(-j)}$

% First, by considering the first-order optimality condition of original optimization problem, an equation about $\hat\beta$ is derived. Analogously, two equations can also be obtained by calculating the first-order optimality condition after leaving out one observation and one predictor respectively. 
  \item [2)]
  Two properties are derived from comparing three version of first-order conditions: a) The $i$-th fitted value $X_i\hat\beta$ has an asymptotic expression composed of two independent random vectors $X_{i,(-j)}$ and $\hat\beta_{(-i),(-j)}$. b) Each coordinate $\hat\beta_j$ can be written as a sum of $n$  random variables which are asymptotically independent. 
%  Second, Subtracting these two equations from original equation and performing Taylor expansion, two approximate equations are derived which characterize the relationship between the estimator after leaving one predictor and one observation respectively.
\item [3)] 
Using above two properties , $\hat\beta_j$ has the same distribution as a combination of several scalar variables when $n,p\to \infty$. Hence every statistic of $\hat\beta$ (such as expectation, variance and first order condition of optimization) can be expressed by these scalar variables, from which the SEs of $\hat\beta$ are derived.
% Then the correlation between the estimator $\hat{\boldsymbol{\beta}}$ and samples $(\rmX, \rvy)$ can be decoupled, which means that the relationship between samples and estimator can be characterized by two independent random vectors. 
% \item [4)] 
%  Hence  previous high-dimensional approximate equations can be easily transformed into 1-dimensional equations when calculating the asymptotic statistics such as MSE and bias as $n,d\to \infty$. Each statistic corresponds to an equation, which can be combined to derive the set of SE.
\end{itemize}

Briefly reviewing the procedures of LOO approach, we find that  the sample matrix $X$ (which is a $\R^{n\times p}$ Gaussian matrix)is decomposed into two independent Gaussian vectors through some special techniques in both LOO and CGMT, which allows the law of large numbers to simplify the first-order equations into scalar equations. This may help us understand the equivalence between CGMT and LOO. The more intrinsic equivalence of these two methods is still under investigation.

% \subsubsection{Equivalence of the reduction form between CGMT and LOO}

%Next we give a heuristic proof of the equivalence between CGMT anf LOO. 
% Briefly reviewing the procedures of LOO approach, we find that the  Gaussian matrix $\rmX$ is decomposed into a form of two independent normal vectors which are derived by leaving one observation and one predictor respectively. A similar decomposition appears in CGMT when simplifying the PO to AO. We believe that this similarity can help us understand the equivalence between the resulting SEs from LOO and CGMT respectively. 

% We notice that $\rmA \vx \sim N(0,\|\vx\|^2 I_m)$, $\rmA \vy \sim N(0,\|\vy\|^2 I_n)$.  Therefore CGMT essentially transforms the originally two correlated random vectors $\rmA\vx, \rmA \vy$ in PO into two independent normal vectors $\|\vx\|\rvg,\|\vy\|\rvh$. Similar to LOO, the solution of AO can be characterized by two independent normal vectors. From this perspective, the CGMT approach has some inherent equivalence with LOO approach. 

\section{Discussion and future directions}\label{Discussion and future directions}
%\subsection{Equivalence of SEs from different methods}\label{Equivalence of SEs from different methods}
In this paper,  we first showed that for the high-dimensional $M$-estimator,  the three sets of SEs derived from AMP, CGMT and LOO are equivalent. We then further show that this equivalence actually appears in various high-dimensional problems. This strongly suggests us that there should be a deep relation between these three approaches. 

Though AMP, CGMT and LOO are different at the first glance, we find that they all can be treated as approximations of the same first order optimality conditions. To be more precise, LOO decouples the correlation between samples and estimator after comparing first-order optimality conditions of the initial  optimization with two leaving-one-out version;
CGMT simplifies the first-order optimality conditions by making some relaxation of the PO problem;
AMP solves the first order optimality conditions directly. All their asymptotic behaviours are characterized by the corresponding SEs respectively.  The equivalence between these SEs sheds us light on looking for a more comprehensive theories to explain this intriguing phenomenon.

\bibliography{reference}
\bibliographystyle{chicago}

\newpage

\vskip .65cm
\noindent
Saidi Luo: First author, Center for Statistical Science, Tsinghua University. 
\vskip 2pt
\noindent
E-mail: lsd18@mails.tsinghua.edu.cn
\vskip 2pt

\noindent
Songtao Tian: Corresponding author, Department of Mathematical Sciences, Tsinghua University.
\vskip 2pt
\noindent
E-mail: tst20@mails.tsinghua.edu.cn
\end{document}

%%%%%%%%%%%%%%%%%%%%%%%%%%%%%%%%%%%%%%%%%%%%%%%%%%%%%%%%%%%%%%%%%%%%%%%%%%%%%%%%%%%%%%%%%%%%%%%%%%%%%%%%%%%%%%%%%%%%%%%%%%%%
%%%%%%%%%%%%%%%%%%%%%%%%%%%%%%%%%%%%%%%%%%%%%%%%%%%%%%%%%%%%%%%%%%%%%%%%%%%%%%%%%%%%%%%%%%%%%%%%%%%%%%%%%%%%%%%%%%%%%%%%%%%%

\renewcommand{\baselinestretch}{2}

\markright{ \hbox{\footnotesize\rm  
}\hfill\\[-12pt]
\hbox{\footnotesize\rm
}\hfill }

\markboth{\hfill{\footnotesize\rm SAIDI LUO, SONGTAO TIAN} \hfill}
{\hfill {\footnotesize\rm Supplementary Material: Equivalence of state equations} \hfill}

\renewcommand{\thefootnote}{}
$\ $\par \fontsize{12}{14pt plus.8pt minus .6pt}\selectfont

%%%%%%%%%%%%%%%%%%%%%%%%%%%%%%%%%%%%%%%%%%%%%%%%%%%%%%%%%%%%%%%%%%%%%%%%%%%%%%%%%%%%%%%%%%%%%%%%%%%%%%%%%%%%%%%%%%%%%%%%%%%%

 \centerline{\large\bf Supplementary Material: Equivalence of state equations from }
\vspace{2pt}
 \centerline{\large\bf different methods in
high-dimensional regression }
\vspace{2pt}
 \author{Author(s)}
\vspace{.4cm}
\renewcommand{\thefootnote}{\fnsymbol{footnote}}
\centerline{ Saidi Luo,  Songtao Tian}
% \footnotetext[1]{The authors contributed equally to this work and are listed in alphabetical order. Saidi Luo: First author, Center for Statistical Science, Tsinghua University. Songtao Tian: Corresponding author, Department of Mathematical Sciences, Tsinghua University.}
\vspace{.4cm} 
\centerline{\it Tsinghua University}
 \vspace{.55cm} \fontsize{9}{11.5pt plus.8pt minus.6pt}\selectfont
\vspace{.55cm}
%  \centerline{\bf Supplementary Material}
\vspace{.55cm}
% \fontsize{9}{11.5pt plus.8pt minus .6pt}\selectfont

\noindent
\large The supplementary material contains the proofs of theorems and propositions.\\
% \par

\setcounter{section}{0}
\setcounter{equation}{0}
\def\theequation{S\arabic{section}.\arabic{equation}}
\def\thesection{S\arabic{section}}

\fontsize{12}{14pt plus.8pt minus .6pt}\selectfont

\section{Proof of Proposition 2.3}
\label{appendix:SEs of M-estimator derived from CGMT}
By the following linear parameter transformation:
$$
\boldsymbol w:=\boldsymbol\beta-\boldsymbol{\beta}^*,
$$
the M-estimator optimization problem  becomes:
\begin{equation}\label{relax-M-former}
\min_{\boldsymbol w}\frac1n\sum_{i=1}^n\rho(\epsilon_i-\rvx_i^T\boldsymbol w).
\end{equation}
Introducing the Lagrange multiplier leads to:
$$
\min_{\boldsymbol w,\boldsymbol v}\max_{\boldsymbol u}\frac{1}{n}\sum_{i=1}^n\rho(v_i)+\frac{1}{\sqrt{n}}u_i(v_i-\epsilon_i+\rvx_i^T\boldsymbol w)
$$
where $\vu:=(u_1,...,u_n);\vv=(v_1,...,v_n)$.
Then we rewrite it in the matrix form:
$$
\min_{\boldsymbol w,\boldsymbol v}\max_{\boldsymbol u}\frac{1}{\sqrt{n}}\boldsymbol u^T\rmX\boldsymbol w+\frac{1}{n}\sum_{i=1}^n\rho(v_i)+\frac{1}{\sqrt{n}}(\boldsymbol u^T\boldsymbol v-\boldsymbol u^T\bm\epsilon)
$$
where $\bm\epsilon:=(\epsilon_1,...,\epsilon_n)^T$. This is just the PO problem in CGMT. 

Denote $\tilde \rmX=\sqrt n\rmX,\tilde \vw=\frac{\vw}{\sqrt n}$, then we have $\tilde X_{i,j}\overset{i.i.d.}{\sim}\mathcal N(0,1),\tilde{\vw}=\frac{\boldsymbol\beta-\boldsymbol\beta^*}{\sqrt n}$. This means $\|\vw\|^2$ is just the MSE of interest and $\tilde \rmX$ is a standard Gaussian matrix composed of iid standard normal variable. However, in the following, we rewrite $\tilde \rmX,\tilde \vw$ as $\rmX, \vw$ respectively for  the simplicity of notation. 

Using CGMT about $\tilde \rmX$ as in \citep{salehi2019impact} by Corollary 3 in it, then the AO problem associated to it is the following min-max problem: 
$$
\min_{\boldsymbol w,\boldsymbol v}\max_{\boldsymbol u}\frac{1}{\sqrt{n}}(||\boldsymbol u||_2\boldsymbol g^T\boldsymbol w+||\boldsymbol w||_2\boldsymbol h^T\boldsymbol u)+\frac{1}{n}\sum_{i=1}^n\rho(v_i)+\frac{1}{\sqrt{n}}(\boldsymbol u^T\boldsymbol v-\boldsymbol u^T\bm\epsilon), 
$$
where $\boldsymbol g \in \R^d$ and $\boldsymbol h \in \R^n$ have i.i.d. $\mathcal{N}(0, 1)$ entries.

%Next, what we want to do is to simplify the above optimization to optimization problem only about serveal scalars. 

Let $\|\boldsymbol w\|_2=\tau_3$, note that now 
\begin{equation}\label{paras stand for MSE in M-estimator from CGMT}
\tau_3^2=\frac{1}{n}||\boldsymbol\beta-\boldsymbol\beta^*||^2 ,
\end{equation}
which is just the MSE.

Then the optimization becomes:
$$
\min_{\tau_3,\boldsymbol v}\max_{\boldsymbol u}\frac{1}{\sqrt{n}}(-\tau_3||\boldsymbol u||_2||\boldsymbol g||_2+\tau_3 \boldsymbol h^T\boldsymbol u)+\frac{1}{n}\sum_{i=1}^n\rho(v_i)+\frac{1}{\sqrt{n}}(\boldsymbol u^T\boldsymbol v-\boldsymbol u^T\bm\epsilon).
$$
Letting $||\boldsymbol u||_2=\mu$, then We have the following optimization:
\begin{equation}\label{relax-M-latter}
\min_{\tau_3,\boldsymbol v}\max_{\mu>0}-\frac{\tau_3\mu}{\sqrt{n}}||\boldsymbol g||_2+\frac{1}{n}\sum_{i=1}^n\rho(v_i)+\frac{\mu}{\sqrt{n}}||\tau \boldsymbol h+\boldsymbol v-\bm\epsilon||_2.
\end{equation}
equivalently:
\begin{equation}\label{equation, equivalently2}
\min_{\tau_3,\boldsymbol v}\max_{\mu>0}\frac{1}{n}\sum_{i=1}^n\rho(v_i)+\mu(\frac{1}{\sqrt{n}}||\tau_3 \boldsymbol h+\boldsymbol v-\bm\epsilon||_2-\frac{\tau_3}{\sqrt{n}}||\boldsymbol g||_2).  
\end{equation}

In order to make the $||\tau_3 \boldsymbol h+\boldsymbol v-\bm\epsilon||_2$ separable, we  use the following optimization:
$$
x=\min_{\alpha>0}\frac{\alpha}{2}+\frac{x^2}{2\alpha}, 
$$
for any $x$ and $\alpha>0$.
Replace $x$ by  $\frac{1}{\sqrt{n}}||\tau_3 \boldsymbol h +\boldsymbol v-\bm\epsilon||_2$, the optimization problem (\ref{equation, equivalently2}) becomes:
$$
\min_{\tau_3,\boldsymbol{v}}\max_\mu-\frac{\tau_3\mu}{\sqrt{n}}||\boldsymbol g||_2+\frac{1}{n}\sum_{i=1}^n\rho(v_i)+\mu(\max_{\alpha>0}\frac{\alpha}{2}+\frac{1}{2\alpha n}||\tau_3 \boldsymbol h+\boldsymbol v-\bm\epsilon||_2^2).
$$

Then what we want to do the scalarization procedure: make the optimization about $\boldsymbol {v}$ becoming optimization about a scalar. First, flipping in the order of min-max by \citep{thrampoulidis2018precise}: 
% and introducing Moreau envelope function $M_\rho(x;t)$ leads to
\begin{equation}\label{equation 25}
\min_{\tau_3,\alpha}\max_\mu\frac{\alpha\mu}{2}-\frac{\tau_3\mu}{\sqrt{n}}||\boldsymbol g||_2+\frac{1}{n}[\min_{\boldsymbol{v}}\sum_{i=1}^n\rho(v_i)+\frac{\mu}{2\alpha}(\tau_3 h_i+v_i-\epsilon_i)^2].    
\end{equation}
Introducing Moreau envelope function $M_\rho(x;t)$, the optimization problem (\ref{equation 25}) becomes:
$$
\min_{\tau_3,\alpha}\max_\mu\frac{\alpha\mu}{2}-\frac{\tau_3\mu}{\sqrt{n}}||\boldsymbol g||_2+\frac{1}{n}\sum_{i=1}^nM_\rho(\epsilon_i-\tau_3 h_i;\alpha/\mu).
$$
By Lemma 9 in Appendix A in \citep{thrampoulidis2018precise}, considering asymptotic $n,p\rightarrow\infty,p/n\rightarrow\kappa^*$ leads to: 
$$
\frac{\alpha\mu}{2}-\frac{\tau_3\mu}{\sqrt{n}}||\boldsymbol g||_2+\frac{1}{n}\sum_{i=1}^n M_\rho(\epsilon_i-\tau_3 h_i;\alpha/\mu)\overset{a.s.}{\rightarrow}\frac{\alpha\mu}{2}-{\tau_3\mu\sqrt{\kappa_*}}+\mathbb EM_\rho(W_3-\tau_3 Z_3;\alpha/\mu)
$$
where $Z_3\sim\mathcal  N(0,1)$ is independent of everything else.

Introduce $W_3\sim P_\epsilon$. Then, asymptotically, we can deal with the following problem:
\begin{equation}\label{equation, SE of M-estimator, min-max problem}
 \min_{\tau_3,\alpha}\max_\mu\frac{\alpha\mu}{2}-{\tau_3\mu}{\sqrt{\kappa_*}}+\mathbb EM_\rho(W_3+\tau_3 Z;\alpha/\mu).  
\end{equation}
Denoting the objective function of (\ref{equation, SE of M-estimator, min-max problem}) by $\phi$, then since $\phi$ is convex about $(\tau_3,\alpha)$ and concave about $\mu$, the saddle point of $\phi$ can be precisely characterized by its first order optimality condition: 
$$
\begin{aligned}
\frac{\partial\phi}{\partial\mu}&=0\Rightarrow\frac{\alpha}{2}-\tau_3\sqrt{\kappa_*}-\frac{\alpha}{\mu^2}\mathbb E[\frac{\partial M_{\rho}}{\partial t}(W_3+\tau_3 Z_3;\alpha/\mu)]=0,\\
\frac{\partial\phi}{\partial\tau}&=0\Rightarrow-\mu\sqrt{\kappa_*}+\mathbb E[Z_3\frac{\partial M_{\rho}}{\partial x}(W_3+\tau_3 Z_3;\alpha/\mu)]=0,\\
\frac{\partial\phi}{\partial\alpha}&=0\Rightarrow\frac{\mu}{2}+\frac{1}{\mu}\mathbb E[\frac{\partial M_{\rho}}{\partial t}(W_3+\tau_3 Z_3;\alpha/\mu)]=0.
\end{aligned}
$$
Combining with (\ref{paras stand for MSE in M-estimator from CGMT}) completes the proof.

\section{Proof of Theorem 1}\label{appendix:Equivalence of SEs of M-estimator from AMP, LOO, CGMT}
 Recall that the SEs of M-estimator from LOO are: 
\begin{equation}
\begin{aligned}
1-\kappa_*&=\mathbb E[\frac{\partial Prox_\rho}{\partial x}(W_1+\tau_1 Z_1;\lambda_1)],\\
\kappa_* \tau_1^2&:=\mathbb E[W_1+\tau_1 Z_1-Prox_\rho(W_1+\tau_1 Z_1;\lambda_1)]^2.
\end{aligned}    
\end{equation}
 SEs of M-estimator from AMP are: 
\begin{equation}
\begin{aligned}
\tau_2^2&=\frac1\kappa_*\lambda_2^2\mathbb E[\frac{\partial M_\rho}{\partial x}(W_2+\tau_2 Z_2;\lambda_2)]^2,\\
\kappa_*&=\lambda_2 \mathbb E[\frac{\partial^2 M_\rho}{\partial x^2}(W_2+\tau_2 Z_2;\lambda_2)].
\end{aligned}    
\end{equation}

Next, we show that these two sets of SEs are equivalent.

First, Simple calculation leads to:
\begin{equation}
\begin{aligned}
\frac{\partial M_\rho(x,t)}{\partial x}&=\rho'(Prox_\rho(x;t)),\\
\frac{\partial^2 M_\rho(x,t)}{\partial x^2}&=\rho''(Prox_\rho(x;t))\frac{\partial Prox_\rho(x;t)}{\partial x}=\frac{\rho''(Prox_\rho(x;t))}{1+t\rho''(Prox_\rho(x;t))},\\
\frac{x-Prox_\rho(x;t)}{t}&=\rho'(Prox_\rho(x;t)).
\end{aligned}
\end{equation}

Combine this, SEs from AMP can be rewritten as:
\begin{equation}
\begin{aligned}
\tau_2^2&=\frac1\kappa_*\lambda_2^2\mathbb E[\rho'(Prox_\rho(W_2+\tau_2 Z_2;\lambda_2))]^2,\\
\kappa_*&=\lambda_2 \mathbb E[\frac{\rho''(Prox_\rho(W_2+\tau_2 Z_2;\lambda_2))}{1+\lambda_2\rho''(Prox_\rho(W_2+\tau_2 Z_2;\lambda_2))}].
\end{aligned}    
\end{equation}
which verify that SEs from LOO and SEs from AMP are equivalent.

Next we prove  equivalence of SEs from AMP and SEs  from CGMT by parameter transformations suggested.

Recall the SEs from CGMT are:
\begin{equation}
\begin{aligned}
\frac{\alpha}{2}-\tau_3\sqrt{\kappa_*}-\frac{\alpha}{\mu^2}\mathbb E[\frac{\partial M_{\rho}}{\partial t}(W_3+\tau_3 Z_3;\alpha/\mu)]=0,\\
-\mu\sqrt{\kappa_*}+\mathbb E[Z_3\frac{\partial M_{\rho}}{\partial x}(W_3+\tau_3 Z_3;\alpha/\mu)]=0,\\
\frac{\mu}{2}+\frac{1}{\mu}\mathbb E[\frac{\partial M_{\rho}}{\partial t}(W_3+\tau_3 Z_3;\alpha/\mu)]=0.
\end{aligned}
\end{equation}
Let $b=\frac{\alpha}{\mu}$, then we have:
\begin{equation}
\begin{aligned}
\frac{\mu b}{2}-\tau_3\sqrt{\kappa_*}-\frac{\alpha}{\mu^2}\mathbb E[\frac{\partial M_{\rho}}{\partial t}(W_3+\tau_3 Z_3;b)]=0,\\
-\mu\sqrt{\kappa_*}+\mathbb E[Z_3\frac{\partial M_{\rho}}{\partial x}(W_3+\tau_3 Z_3;b)]=0,\\
\frac{\mu}{2}+\frac{1}{\mu}\mathbb E[\frac{\partial M_{\rho}}{\partial y}(W+\tau Z;b)]=0.
\end{aligned}
\end{equation}

Comparing these results leads to:
\begin{equation}
   \mathbb E[\frac{\partial M_{\rho}}{\partial x}(W_3+\tau_3 Z_3;b)]^2=\frac{\tau_3^2\kappa_*}{b^2}\\
\mathbb E[\frac{\partial^2 M_{\rho}}{\partial x^2}(W_3+\tau_3 Z_3;b)]=\frac{\kappa_*}{b} .
\end{equation}

Combining stein lemma and
\begin{equation}
\frac{\partial M_\rho}{\partial t}(W_3+\tau_3 Z_3;b)=-\frac{1}{2}[\frac{\partial M_\rho}{\partial x}(W_3+\tau_3 Z_3;b)]^2
\end{equation}
completes our proof.
\section{Proof of Proposition 3.2}\label{appendix:Equivalence of SEs of LASSO from AMP and CGMT}
%You may include other additional sections here.

In this proof, we refer to the technique developed in \citep{donoho2016high}.  
The LASSO problem solves
% Since only the solution of optimization is interested, we first transform LASSO optimization under the condition without changing the solution.
% The solution is not changed when an optimization is divided by $n$. Hence the LASSO problem solves
\begin{equation}
\arg\min_{\boldsymbol{\beta}} \frac{\lambda_*}{n}\|\boldsymbol{\beta}\|_1+ \frac{1}{2n}\|\mathbf{y}-\mathbf{X}\boldsymbol{\beta}\|^2.
\end{equation}

Notice that $\|\mathbf{y}-\mathbf{X}\boldsymbol{\beta}\|^2 = \sum_{i=1}^n\left[(x_i^T\boldsymbol{\beta})^2-2y_ix_i^T\boldsymbol{\beta}+ y_i^2\right]$ .
The optimization can be transformed into
% and noticing that $\sum_{i=1}^n y_i^2$ is unchanged with respect to $\boldsymbol{\beta}$, we can turn the optimization into 
\begin{equation}
\arg\min_{\boldsymbol{\beta}} \frac{\lambda_*}{n}\|\boldsymbol{\beta}\|_1+ \frac{1}{n}\sum_{i=1}^n\left[ \frac{1}{2}(x_i^T\boldsymbol{\beta})^2-y_i x_i^T\boldsymbol{\beta} \right].
\end{equation}

In the following proof, we consider a more general optimization than LASSO:
\begin{equation}
\arg\min_{\boldsymbol{\beta}} \frac{\lambda_*}{n}f(\boldsymbol{\beta})+ \frac{1}{n}\sum_{i=1}^n\left[ \mathbf{1}^T\bm{\rho}(\mathbf{u})-y_i x_i^T\boldsymbol{\beta} \right].
\end{equation}
where $\bm{\rho}(\cdot)$ and $f(\cdot)$ are 'separable' in the sense that there exist  scalar functions $\rho(\cdot), \tilde f(\cdot)$ so that $\bm{\rho}(\cdot)$ and $f(\cdot)$ can be expressed as the following form:   $\boldsymbol{\rho}(\mathbf{x}) = (\rho(x_1),\cdots,\rho(x_d))^T$ and $f(\mathbf{x}) = \sum_{i=1}^d \tilde f(x_i)$. In particular, in LASSO, $\rho(t) = \frac 12 t^2$ and $f(\mathbf{x}) = \|\mathbf{x}\|_1$. 

In order to apply CGMT, we introduce a new variable $\vu$ and have following optimization
\begin{equation} 
\begin{aligned}
\min_{\boldsymbol{\beta},\mathbf{u}}& \frac{\lambda_*}{n}f(\boldsymbol{\beta})+\frac{1}{n}\left(\mathbf{1}^T\rho(\mathbf{u})-y^T\mathbf{u}\right),\\
\text{s.t.}\quad& \mathbf{u}=\mathbf{X}\boldsymbol{\beta}=\frac{1}{\sqrt{n}}\mathbf{H}^*\boldsymbol{\beta}, \label{AppdixLasso}
\end{aligned} 
\end{equation} 
where $\mathbf{H}^*=\sqrt{n}\mathbf{X}\in \mathbb{R}^{n\times d}$ and hence $\mathbf{H}^*_{ij} \overset{\text{i.i.d. }}{\sim} \mathcal{N}(0,1)$. By using Lagrange multiplier we can rewrite (\ref{AppdixLasso}) as a min-max optimization: 
\begin{equation} \label{17}
\min_{\boldsymbol{\beta}\in\R^d,\mathbf{u}\in\mathbb{R}^n}\max_{\vv} \frac 1n \mathbf{1}^T \rho(\mathbf{u})-\frac{1}{n}y^T\mathbf{u} + \frac{\lambda_*}{n}\|\boldsymbol{\beta}\|_1+\frac{1}{n}\vv^T(\mathbf{u}-\frac{1}{\sqrt{n}}\mathbf{H^*}\boldsymbol{\beta}).
\end{equation}

Denote $P = \frac{\boldsymbol{\beta}^* (\boldsymbol{\beta}^*)^T}{\|\boldsymbol{\beta}^*\|^2}$ as the projection matrix of true signal $\boldsymbol{\beta}^*$ and $P^\perp = I_d - P$ as the orthogonal complement. To apply CGMT, we need first decompose $\mathbf{H^*}$ into
\begin{equation}
\begin{aligned}
\mathbf{H}_1^*=\mathbf{H}^*\cdot P,&\quad \mathbf{H}_2^*=\mathbf{H}^*\cdot P^\perp ,\\
\mathbf{H}^*=&\mathbf{H}_1^*+\mathbf{H}_2^*.
\end{aligned}
\end{equation}
In addition, Recalling the linear model:
\begin{equation}
\begin{aligned}
\ry_i = \rvx_i^T \boldsymbol{\beta}^* +\epsilon_i, \quad \text{for } i\in [n]  .  
\end{aligned}
\end{equation}
we have $y=\mathbf{X}\boldsymbol{\beta}^*+\bm\epsilon = \frac{1}{\sqrt{n}}\mathbf{H}^*\boldsymbol{\beta}^* = \frac{1}{\sqrt{n}}\mathbf{H}_1^*\boldsymbol{\beta}^*$. Hence (\ref{17}) can be rewritten as
\begin{equation}
\min_{\boldsymbol{\beta}\in\R^d,\mathbf{u}\in\mathbb{R}^n}\max_{\vv} \frac{1}{n}\mathbf{1}^T\rho(\mathbf{u})-\frac{1}{n}y^T\mathbf{u} + \frac{\lambda_*}{n}f(\boldsymbol{\beta})+\frac{1}{n}\vv^T(\mathbf{u}-\frac{1}{\sqrt{n}}\mathbf{H}_1^*\boldsymbol{\beta}) - \frac{1}{n\sqrt{n}}\vv^T\mathbf{H}_2^*\boldsymbol{\beta}. \label{41}
\end{equation}
By using CGMT for $\vv^T\mathbf{H}_2^*\boldsymbol{\beta}$ as in \citep{salehi2019impact} by Corollary 3 in it, the corresponding AO of (\ref{41}) is
\begin{equation}
\begin{aligned}
\min_{\boldsymbol{\beta}\in\R^d,\mathbf{u}\in\mathbb{R}^n}\max_{\vv} &\frac{1}{n}\mathbf{1}^T\rho(\mathbf{u})-\frac{1}{n}y^T\mathbf{u} + \frac{\lambda_*}{n}f(\boldsymbol{\beta})+\frac{1}{n}\vv^T(\mathbf{u}-\frac{1}{\sqrt{n}}\mathbf{H}_1^*\boldsymbol{\beta})\\
&-\frac{1}{n\sqrt{n}}\left(\vv^Th\|P^\perp\boldsymbol{\beta}\|_2+\|\vv\|_2g^TP^\perp\boldsymbol{\beta}\right).
\end{aligned}
\label{AO}
\end{equation}
where $h\sim \mathcal{N}(0,I_n)$ and $g\sim \mathcal{N}(0,I_d)$ are two independent gaussian vectors.

We first consider the maximization with respect to the direction of $\vv$.  The part related to $\vv$ in optimization (\ref{AO}) is:
\begin{equation}
\max_{\vv\in\R^n} \frac{1}{n\sqrt{n}}\|\vv\|_2g^TP^\perp\boldsymbol{\beta}+\frac{1}{n}\vv^T\left(\mathbf{u}-\frac{1}{\sqrt{n}}\mathbf{H}_1^*\boldsymbol{\beta}-\frac{1}{\sqrt{n}}h\|P^\perp\boldsymbol{\beta}\|_2 \right). \label{44}
\end{equation}

Denoting $r:= \frac {\|\vv\|_2}{\sqrt n}$ and maximizing along the direction of $\vv$ give  
\begin{equation}
\max_{r\ge0} r\left(\frac{1}{n}g^TP^\perp\boldsymbol{\beta}+\|\frac{1}{\sqrt{n}}\mathbf{u}-\frac{1}{n}\mathbf{H}_1^*\boldsymbol{\beta}-\frac{\|P^\perp\boldsymbol{\beta}\|_2}{n}h \|_2 \right) \label{45}.
\end{equation}

% Replacing this in (\ref{AO}) gives

Inserting this into  (\ref{AO}) gives
\begin{equation}
\min_{\boldsymbol{\beta},\mathbf{u}}\max_{r\ge0} \frac{1}{n}\mathbf{1}^T\rho(\mathbf{u})-\frac{1}{n}y^T\mathbf{u} + \frac{\lambda_*}{n}f(\boldsymbol{\beta}) 
+r\left(\frac{1}{n} g^TP^\perp\boldsymbol{\beta}
+\left\| \frac{1}{\sqrt{n}}\mathbf{u}-\frac{1}{n}\mathbf{H}_1^*\boldsymbol{\beta}-\frac{\|P^\perp\boldsymbol{\beta}\|_2}{n}h \right\|_2 \right). \label{46}
\end{equation}

% In addition, we introduce new variables $\boldsymbol{\mu}, \vw$ in which $\boldsymbol{\mu}$ replace $\boldsymbol{\beta}$ in $\|\boldsymbol{\beta}\|_1$ and $\vw$ viewed as Lagrangian. 
In addition, introduce $\boldsymbol{\mu}$ to replace $\boldsymbol{\beta}$ in $\|\boldsymbol{\beta}\|_1$ Lagrangian and $\vw$. 
Then  (\ref{46}) can be rewritten as, 
\begin{equation}
\begin{aligned}
\min_{\mathbf{u}\in\R^n; \atop \boldsymbol{\beta},\boldsymbol{\mu}\in\R^d} \max_{r\ge0; \atop \vw\in\R^d} 
&\frac{1}{n}\mathbf{1}^T\rho(\mathbf{u})-\frac{1}{n}y^T\mathbf{u} + \frac{\lambda_*}{n}f(\boldsymbol{\mu}) + r\left(\frac{1}{n}g^TP^\perp\boldsymbol{\beta} \right)
\\
 &+r\left\| \frac{1}{\sqrt{n}}\mathbf{u}-\frac{1}{n}\mathbf{H}_1^*\boldsymbol{\beta}-\frac{\|P^\perp\boldsymbol{\beta}\|_2}{n}h \right\| _2
 +\frac{1}{d} \vw^T(\boldsymbol{\mu}-\boldsymbol{\beta}).
 \end{aligned} \label{47}
 \end{equation}

We define $\alpha = \frac{\boldsymbol{\beta}^T\boldsymbol{\beta}^*}{\|\boldsymbol{\beta}^*\|_2^2}$,  $\sigma = \frac{1}{\sqrt{n}}\|P^\perp \boldsymbol{\beta}\|_2$ and $\mathbf{q} = \frac{\mathbf{H^*\boldsymbol{\beta}^*}}{r_{1*}\sqrt{n}}$ where $r_{1*} = \frac{\|\boldsymbol{\beta}^*\|_2}{\sqrt{n}}$. Then $\mathbf{q}$ is a standard Gaussian vector and 
\begin{equation} \label{48}
\frac{1}{n}\mathbf{H}_1^*\boldsymbol{\beta} = \frac{1}{n}\mathbf{H}^*(P\boldsymbol{\beta}) = \frac{\mathbf{H}^*}{n}\cdot\alpha\boldsymbol{\beta}^*\overset{d}=\frac{\alpha}{n}r_{1*}\sqrt{n}\mathbf{q}.
\end{equation} 
Decomposing $\vw = (P+P^{\perp})\vw$, then the last item in (\ref{47}) can be rewritten as
\begin{equation}\label{equation, decomposition of w^T(mu-beta)}
\frac{1}{d}\vw^T(\boldsymbol{\mu}-\boldsymbol{\beta}) = \frac{1}{d}(P\vw)^T\boldsymbol{\mu}+\frac{1}{d}(P^\perp \vw)^T\boldsymbol{\mu} 
-\frac{1}{d}(P\vw)^T\boldsymbol{\beta} - \frac{1}{d}(P^\perp \vw)^T\boldsymbol{\beta}.
\end{equation}
Inserting (\ref{48}) and (\ref{equation, decomposition of w^T(mu-beta)}) into (\ref{47}), we have,
\begin{equation} \label{47.1}
\begin{aligned}
\min_{\mathbf{u}\in\R^n; \atop \boldsymbol{\beta},\boldsymbol{\mu}\in\R^d} \max_{r\ge0; \atop \vw\in\R^d} 
&\frac{1}{n}\mathbf{1}^T\rho(\mathbf{u})-\frac{1}{n}y^T\mathbf{u} + \frac{\lambda_*}{n}f(\boldsymbol{\mu}) + r\left(\frac{1}{n}g^TP^\perp\boldsymbol{\beta} \right)- \frac{1}{d}(P^\perp \vw)^T\boldsymbol{\beta}
\\
 &+r\left\| \frac{1}{\sqrt{n}}\mathbf{u}-\frac{\alpha r_{1*}}{\sqrt{n}}\mathbf{q}-\frac{\sigma}{\sqrt{n}}h \right\|_2 
 +\frac{1}{d}(P\vw)^T\boldsymbol{\mu}+\frac{1}{d}(P^\perp \vw)^T\boldsymbol{\mu} 
-\frac{1}{d}(P\vw)^T\boldsymbol{\beta}.
 \end{aligned} 
 \end{equation}

Then we can fix $P\boldsymbol{\beta}$ and consider the minimization along the direction of $P^\perp \boldsymbol{\beta}$. Considering the optimization  related to $P^\perp \boldsymbol{\beta}$, we have  
\begin{equation}\label{equation, min Pverticalbeta}
\begin{aligned}
\min_{P^\perp\boldsymbol{\beta}}\frac{r}{n}g^TP^\perp\boldsymbol{\beta}-\frac{1}{d}(P^\perp \vw)^T\boldsymbol{\beta} 
&= \min_{P^\perp\boldsymbol{\beta}}(\frac{r}{n}g^T-\frac{1}{d}\vw^T) P^\perp\boldsymbol{\beta} \\
&= -\|P^\perp\boldsymbol{\beta}\|_2 \cdot\|\frac{r}{n}P^\perp g-\frac{1}{d}P^\perp \vw \|_2 \\
&= -\sigma \cdot \|\frac{r}{\sqrt{n}}P^\perp g-\sqrt{\frac{1}{d\kappa_*}}P^\perp \vw\|_2.
\end{aligned}
\end{equation}

Notice that (\ref{47.1}) reaches optimal when $\boldsymbol{\mu} = \boldsymbol{\beta}$. Inserting (\ref{equation, min Pverticalbeta}) into (\ref{47.1}) leads to
\begin{equation}
\begin{aligned}
\min_{\mathbf{u}\in\R^n,\boldsymbol{\mu}\in\R^d; \atop  \alpha\in\R,\sigma\ge0 } \max_{r\ge0; \atop \vw\in\R^d} 
&\frac{1}{n}\mathbf{1}^T\rho(\mathbf{u})-\frac{1}{n}y^T\mathbf{u} + \frac{\lambda_*}{n}f(\boldsymbol{\mu}) 
-\sigma \cdot \left\|\frac{r}{\sqrt{n}}P^\perp g-\sqrt{\frac{1}{d\kappa_*}}P^\perp \vw\right\|
\\
&+r\left\| \frac{1}{\sqrt{n}}\mathbf{u}-\frac{\alpha r_{1*}}{\sqrt{n}}\mathbf{q}-\frac{\sigma}{\sqrt{n}}h \right\| 
+\frac{1}{d}(P^\perp \vw)^T\boldsymbol{\mu}.
\end{aligned} \label{47.2}
\end{equation}

For the  simplifying procedures in the following steps in our analysis, we change  $\|\cdot\|_2 \to \|\cdot\|_2^2$ by
\begin{equation} \label{30}
\begin{aligned}
rx &= \min_{v\ge0} \frac{r}{2v}+\frac{rv}{2}x^2, \\
-\sigma x &= \max_{\tau\ge0} -\frac{\sigma}{2\tau}-\frac{\sigma\tau}{2}x^2.
\end{aligned}
\end{equation}

Applying (\ref{30}), we are able to rewrite (\ref{47.2}) as
\begin{equation}
\begin{aligned}
\min_{\mathbf{u}\in\R^n,\boldsymbol{\mu}\in\R^d; \atop
\alpha\in\R,v,\sigma\ge0 } \max_{\vw\in\R^d;\atop r,\tau\ge0} 
&\frac{1}{n}\mathbf{1}^T\rho(\mathbf{u})-\frac{1}{n}y^T\mathbf{u} + \frac{\lambda_*}{n}f(\boldsymbol{\mu}) 
-\frac{\sigma}{2\tau} -\frac{\sigma\tau}{2} \|\frac{r}{\sqrt{n}}P^\perp g-\sqrt{\frac{1}{d\kappa_*}}P^\perp \vw\|^2
\\
&+\frac{r}{2v} + \frac{rv}{2}\left\| \frac{1}{\sqrt{n}}\mathbf{u}-\frac{\alpha r_{1*}}{\sqrt{n}}\mathbf{q}-\frac{\sigma}{\sqrt{n}}h \right\|^2 
+\frac{1}{d}(P^\perp \vw)^T\boldsymbol{\mu}.
\end{aligned} \label{49}
\end{equation}

{\bf Optimization with respect to $\vw$}:  Next we consider the maximization with respect to $\vw$. We first extract the item related to $\vw$ in (\ref{49}) and apply the completion of squares:
\begin{equation}\label{32}
\begin{aligned}
\max_\vw& -\frac{\sigma\tau}{2}\left\|\frac{r}{\sqrt{n}}P^\perp g-\sqrt{\frac{1}{d\kappa_*}}P^\perp \vw  \right\|^2 
+ \frac{1}{d}(P^\perp \vw)^T\boldsymbol{\mu}  \\
= \max_{\vw}& -\frac{\sigma\tau}{2}\left\|\frac{r}{\sqrt{n}}P^\perp g-\sqrt{\frac{1}{d\kappa_*}}P^\perp \vw +\frac{1}{\sqrt{d/\kappa_*}\sigma\tau}P^\perp \boldsymbol{\mu} \right\|^2 \\
&+\frac{1}{2n\sigma\tau}\left\|P^\perp\boldsymbol{\mu}+\sigma\tau rP^\perp g \right\|^2 
-\frac{\sigma\tau r^2}{2n}\|P^\perp g\|^2.
\end{aligned} 
\end{equation} 

\begin{itemize}
\item[1)] For the last item in (\ref{32}), since $g\sim \mathcal{N}(0,I_d)$ and $P^\perp$ is a $(n-1)$-dimensional projection matrix,  we derive that $\|P^\perp g\|_2^2 \sim \left\|d(0,(P^\perp)^2)\right\|_2^2 \overset{d}{=}\chi_{d-1}^2$ and 
\begin{equation}
\frac{\sigma\tau r^2}{2n}\|P^\perp g\|^2 \stackrel{a.s.}{\to} 
\frac{\sigma\tau r^2\kappa_*}{2}.
\end{equation}

\item[2)]Since $P^\perp = I_d-P$, the second item in (\ref{32}) can be rewritten as
\begin{equation} \label{34}
\begin{aligned}
\frac{1}{n}\left\|P^\perp\boldsymbol{\mu}+\sigma\tau rP^\perp g \right\|^2 =& \frac1n\|\boldsymbol{\mu}+\sigma\tau rg\|^2-\frac1n\|P\boldsymbol{\mu}\|^2\\
&-\frac{(\sigma\tau r)^2}{n}\|Pg\|^2  - \frac{2\sigma\tau r}{n}
(Pg)^T\boldsymbol{\mu}.
\end{aligned} 
\end{equation}

The last two items of (\ref{34}) can be omitted in the limit of $d,n \to \infty$ because $\frac{\|Pg\|^2}{n} = O_p(\frac 1n)$ and $\frac{1}{n}(Pg)^T\boldsymbol{\mu} = O_p(\frac 1{\sqrt n})$. 
The second item of (\ref{34}) is $\frac 1n\|P\boldsymbol{\mu}\|^2 = \frac 1n\|P\boldsymbol{\beta}\|^2 = \alpha^2{r_{1*}}^2$ by definition.
\item[3)]The first item in (\ref{32}) reaches 0 when maximizing $\vw$.
\end{itemize}
The optimization (\ref{49}) now can be rewritten as:
\begin{equation} \label{51}
\begin{aligned}
\min_{\mathbf{u}\in\R^n,\boldsymbol{\mu}\in\R^d;\atop \alpha\in\R,v,\sigma\ge0 } \max_{r,\tau\ge0} 
&\frac{1}{n}\mathbf{1}^T\rho(\mathbf{u})-\frac{1}{n}y^T\mathbf{u} + \frac{\lambda_*}{n}f(\boldsymbol{\mu}) 
-\frac{\sigma}{2\tau} 
+\frac{1}{2n\sigma\tau}\|\boldsymbol{\mu}+\sigma\tau rg\|^2-\frac{\alpha^2r_{1*}^2}{2\sigma\tau} 
- \frac{\sigma\tau r^2\kappa_*}{2}
\\
&+\frac{r}{2v} + \frac{rv}{2}\left\| \frac{1}{\sqrt{n}}\mathbf{u}-\frac{\alpha r_{1*}}{\sqrt{n}}\mathbf{q}-\frac{\sigma}{\sqrt{n}}h \right\|^2. 
\end{aligned} 
\end{equation}

{\bf Optimization respect to $\boldsymbol{\mu}$}: Consider the items related to $\boldsymbol{\mu}$ in (\ref{51}) :
\begin{equation} \label{52}
\begin{aligned}
\min_{\boldsymbol{\mu}\in\R^d} \frac{\lambda_*}{n}f(\boldsymbol{\mu})&+\frac{1}{2n\sigma\tau}\|\boldsymbol{\mu}+\sigma\tau rg\|_2^2 ,\\
s.t.& \qquad\ \boldsymbol{\mu} = \boldsymbol{\beta}.
\end{aligned}
\end{equation}

Notice that $g\sim \mathcal{N}(0,I_d)$,  $\|\boldsymbol{\mu}+\sigma\tau rg\|_2^2\overset{d}{=} \|\boldsymbol{\mu}-\sigma\tau rg\|_2^2$. We rewrite (\ref{52}) as:
\begin{equation} \label{53}
\begin{aligned}
\min_{\boldsymbol{\mu}\in\R^d} \frac{\lambda_*}{n}f(\boldsymbol{\mu})
+\frac{1}{2n\sigma\tau}\|\boldsymbol{\mu}-\sigma\tau rg\|_2^2, \\
s.t. \qquad\ 
\frac{1}{n}{\boldsymbol{\beta}^*}^T\boldsymbol{\mu} = \frac{1}{n}{\boldsymbol{\beta}^*}^T\boldsymbol{\beta} 
= \frac{1}{n}n\alpha r_{1*}^2 = \alpha r_{1*}^2.
\end{aligned}
\end{equation} 

Introducing Lagrangian $\theta$, (\ref{53}) can be rewrite as:
\begin{equation} \label{54}
\min_{\boldsymbol{\mu}\in \R^d}\max_{\theta\in\R} 
\frac{1}{2n\sigma\tau}\|\boldsymbol{\mu} - \sigma\tau rg\|^2+\frac{\lambda_*}{n}f(\boldsymbol{\mu})
- \frac{\theta}{n}{\boldsymbol{\beta}^*}^T\boldsymbol{\mu} + \alpha\theta r_{1*}^2. 
\end{equation} 

Applying the completion of squares to 1,3 items in (\ref{54}) we have:
\begin{equation} \label{54.5}
\begin{aligned}
\frac{1}{2n\sigma\tau}\|\boldsymbol{\mu}-\sigma\tau rg\|^2
- \frac{\theta}{n}{\boldsymbol{\beta}^*}^T\boldsymbol{\mu} =& \frac{1}{2n\sigma\tau}\|\boldsymbol{\mu}-\sigma\tau rg-\theta\sigma\tau\boldsymbol{\beta}^* \|^2 \\
& - \frac{(\theta\sigma\tau)^2}{2n\sigma\tau}\|\boldsymbol{\beta}^*\|^2
-\frac{\theta r\sigma^2\tau^2}{2n\sigma\tau} g^T\boldsymbol{\beta}^*.
\end{aligned}
\end{equation}

The third item can be omitted since $\frac{g^T\boldsymbol{\beta}^*}{n}  = O_p(\frac 1{\sqrt n})$ and the second item has limit
$- \frac{(\theta\sigma\tau)^2}{2n\sigma\tau}\|\boldsymbol{\beta}^*\|^2 \to - \frac{(\theta\sigma\tau)^2}{2n\sigma\tau}\cdot nr_{1*}^2 = -\frac{\sigma\tau\theta^2r_{1*}^2}{2}$. Hence we rewrite right side of (\ref{54.5}) as
\begin{equation}
\frac{1}{2n\sigma\tau}\|\boldsymbol{\mu}-\sigma\tau rg\|^2
- \frac{\theta}{n}{\boldsymbol{\beta}^*}^T\boldsymbol{\mu} = \frac{1}{2n\sigma\tau}\|\boldsymbol{\mu}-\sigma\tau rg-\theta\sigma\tau\boldsymbol{\beta}^* \|^2-\frac{\sigma\tau\theta^2r_{1*}^2}{2}. \label{55}
\end{equation}

Next, denote $\tilde f(x)$ as single-entry form of $f(x)$. We can rewrite the (\ref{55}) in terms of Moreau envelope entry-wisely as follows
\begin{equation}
\begin{aligned}
&\min_{\boldsymbol{\mu}\in \R^d}\max_{\theta\in\R} 
\frac{1}{2n\sigma\tau}\|\boldsymbol{\mu} - \sigma\tau rg\|^2+\frac{\lambda_*}{n}f(\boldsymbol{\mu})
- \frac{\theta}{n}{\boldsymbol{\beta}^*}^T\boldsymbol{\mu} + \alpha\theta r_{1*}^2 \\
=& \max_{\theta} \frac1n M_{\lambda_* \tilde{f}}(\sigma\tau(rg+\theta\boldsymbol{\beta}^*);\sigma\tau)
+ \alpha\theta r_{1*}^2 - \frac{\sigma\tau\theta^2 r_{1*}^2}{2}.
\end{aligned} \label{56.2}
\end{equation}

Substituting (\ref{56.2}) in (\ref{51}) we have,
\begin{equation}
\begin{aligned}
\min_{\mathbf{u}\in\R^n; \atop \alpha\in\R,\sigma,v\ge0} \max_{r,\tau\ge0; \atop \theta\in \R} 
& \frac{1}{n}\mathbf{1}^T\rho(\mathbf{u})-\frac{1}{n}y^T\mathbf{u} 
+ \frac{rv}{2}\left\| \frac{1}{\sqrt{n}}\mathbf{u}-\frac{\alpha r_{1*}}{\sqrt{n}}\mathbf{q}-\frac{\sigma}{\sqrt{n}}h \right\|^2
-\frac{\sigma}{2\tau} -\frac{\alpha^2r_{1*}^2}{2\sigma\tau} \\
&- \frac{\sigma\tau r^2\kappa_*}{2} + \frac{r}{2v} 
+\frac{1}n M_{\lambda_* f}(\sigma\tau(rg+\theta\boldsymbol{\beta}^*);\sigma\tau)
+\alpha\theta r_{1*}^2 - \frac{\sigma\tau\theta^2r_{1*}^2}{2}. 
\end{aligned} \label{57}    
\end{equation}

{\bf Optimization respect to $\vu$}: First we consider the items related to $\vu$. The optimization is 
\begin{equation}
\min_{\mathbf{u}\in\R^n} 
\frac{1}{n}\mathbf{1}^T\rho(\mathbf{u}) - \frac{1}{n}y^T\mathbf{u}
+ \frac{rv}{2}\left\| \frac{1}{\sqrt{n}}\mathbf{u}-\frac{\alpha r_{1*}}{\sqrt{n}}\mathbf{q}-\frac{\sigma}{\sqrt{n}}h \right\|^2. \label{58}
\end{equation}

Applying the completion of squares we have,
\begin{equation}
\begin{aligned}
-\frac{1}{n} y^{T} \mathbf{u}
+\frac{r v}{2}
\left\|\frac{1}{\sqrt{n}} \mathbf{u}-\frac{\alpha r_{1*} }{\sqrt{n}} \mathbf{q}-\frac{\sigma}{\sqrt{n}} h\right\|^{2} 
=& \frac{r v}{2}\left\|\frac{1}{\sqrt{n}} \mathbf{u}-\frac{\alpha r_{1*} }{\sqrt{n}} \mathbf{q}-\frac{\sigma}{\sqrt{n}} h-\frac{1}{r v \sqrt{n}} y\right\|^{2} \\
&-\frac{1}{2rvn}\|y\|^{2}-\frac{r_{1*} \alpha}{n} y^{T} \mathbf{q}-\frac{\sigma}{n} y^{T} h.
\end{aligned} \label{59}
\end{equation}

Using the strong law of large numbers we have,
\begin{equation}\label{60}
\begin{aligned}
    -\frac1{2rvn}\|y\|^2&\stackrel{a.s.}{\to}-\frac{r_{1*}^2+\sigma_*^2}{2rv},\\
    -\frac{r_{1*} \alpha}{n} y^{T} \mathbf{q} &\stackrel{a.s.}{\to} -r_{1*}^2\alpha,\\
    -\frac\sigma n y^Th &\stackrel{a.s.}{\to} 0. 
\end{aligned}
\end{equation}

Next, by substituting (\ref{59}), (\ref{60}) in (\ref{58}), we can rewritten the optimization as, 
\begin{equation}
\begin{aligned}
\min_{\mathbf{u}\in\R^n} 
&\frac{1}{n}\mathbf{1}^T\rho(\mathbf{u}) - \frac{1}{n}y^T\mathbf{u}
+ \frac{rv}{2}\left\| \frac{1}{\sqrt{n}}\mathbf{u}-\frac{\alpha r_{1*}}{\sqrt{n}}\mathbf{q}-\frac{\sigma}{\sqrt{n}}h \right\|^2 \\
= \min_{\mathbf{u}\in\R^n} 
&\frac{1}{n}\mathbf{1}^T\rho(\mathbf{u}) 
+ \frac{r v}{2}\left\|\frac{1}{\sqrt{n}} \mathbf{u}-\frac{\alpha r_{1*} }{\sqrt{n}} \mathbf{q}-\frac{\sigma}{\sqrt{n}} h-\frac{1}{r v \sqrt{n}} y\right\|^{2} \\
 &-\frac{r_{1*}^2+\sigma_*^2}{2rv} - r_{1*}^2\alpha.
\end{aligned} \label{62}
\end{equation}

Then we can rewrite (\ref{62}) in terms of Moreau envelope,
\begin{equation}
\begin{aligned}
&\min_{\mathbf{u}\in \R^n}
\frac{1}{n}\mathbf{1}^T\rho(\mathbf{u}) 
+ \frac{r v}{2}\left\|\frac{1}{\sqrt{n}} \mathbf{u}-\frac{\alpha r_{1*} }{\sqrt{n}} \mathbf{q}-\frac{\sigma}{\sqrt{n}} h-\frac{1}{r v \sqrt{n}} y\right\|^{2} \\
=&\frac1n M_{\rho(\cdot)}(\alpha r_{1*}\mathbf{q}+\sigma h+\frac1{rv}y;\frac1{rv}).
\end{aligned}
\label{62.1}
\end{equation}

Substituting (\ref{62}), (\ref{62.1}) in (\ref{57}) we have
\begin{equation}
\begin{aligned}
\min_{\alpha\in\R; \atop \sigma,v\ge0}
\max_{r,\tau\ge0;\atop \theta\in \R}
\quad &\frac1n M_{\rho}(\alpha r_{1*}\mathbf{q}+\sigma h+\frac1{rv}y;\frac1{rv}) 
-\frac{r_{1*}^2+\sigma_*^2}{2rv} - r_{1*}^2\alpha \\
-&\frac{\sigma}{2\tau} -\frac{\alpha^2r_{1*}^2}{2\sigma\tau} 
- \frac{\sigma\tau r^2\kappa_*}{2} + \frac{r}{2v}
+\alpha\theta r_{1*}^2 - \frac{\sigma\tau\theta^2 r_{1*}^2}{2} \\
+&\frac{1}n M_{\lambda_* f}(\sigma\tau(rg+\theta\boldsymbol{\beta}^*);\sigma\tau).
\end{aligned} \label{63}
\end{equation}

{\bf Final scalarization}: Using the strong law of large number ($\mathbf{q}, h,y,g,\boldsymbol{\beta}^*$ are entry-wise i.i.d.), we can rewrite (\ref{63}) as
\begin{equation}
\begin{aligned}
\min_{\alpha\in\R, \atop \sigma,v\ge0}
\max_{r,\tau\ge0; \atop \theta\in \R}
\quad &-\frac{\sigma}{2\tau} -\frac{\alpha^2r_{1*}^2}{2\sigma\tau} 
- \frac{\sigma\tau r^2 \kappa_*}{2} + \frac{r}{2v}
+\alpha\theta r_{1*}^2 - \frac{\sigma\tau\theta^2r_{1*}^2}{2} -\frac{r_{1*}^2+\sigma_*^2}{2rv} - r_{1*}^2\alpha\\
&+ \mathbb{E}\left[ M_{\rho}(\alpha r_{1*}Z_1+\sigma Z_2+\frac1{rv}(r_{1*}Z_1+\sigma_*Z_3);\frac1{rv})\right] \\
&+ \mathbb{E}\left[ M_{\lambda_* \tilde{f}}(\sigma\tau(rZ+\theta b_0);\sigma\tau)\right]\cdot\frac dn,
\end{aligned}
 \label{67}
\end{equation}
where $Z_1,Z_2,Z \sim \mathcal{N}(0,1)$, $\sigma_* Z_3 \sim P_{\epsilon}$ and  $b_0 \sim \Pi$ are all independent.  

For LASSO, $\tilde f(x) = |x|$. The Moreau envelope $M_{\lambda\tilde{f}}(\cdot;\cdot)$ has property:
\begin{equation}
M_{\lambda_*\tilde{f}}(\sigma\tau(rZ+\theta b_0);\sigma\tau) = \sigma\tau\cdot M_{\lambda_*\tilde{f}}(rZ+\theta b_0;1).
\end{equation}
Besides, for $\rho(x) = \frac 12 x^2$ in LASSO, the Moreau envelope $M_\rho(\cdot;\cdot)$ has explict form:
\begin{equation}
 M_{\rho}(v;t) = \frac{v^2}{2(t+1)},    
\end{equation}
and the second last item of (\ref{67}) can be simplified to:
\begin{equation}
\begin{aligned}
&\mathbb{E}\left[ M_{\rho(\cdot)}(\alpha r_{1*}Z_1+\sigma Z_2+\frac1{rv}(r_{1*}Z_1+\sigma_*Z_3);\frac1{rv})\right]\\ =& \mathbb{E}\left[
 \frac{\left((\alpha+\frac1{rv})r_{1*}Z_1 + \sigma Z_2+\frac{\sigma_*}{rv}Z_3 \right)^2}{2(\frac1{rv}+1)} \right]\\
=&\frac{r_{1*}^2(\alpha rv+1)^2+r^2v^2\sigma^2+\sigma_*^2}{2(1+rv)rv}.
\end{aligned} \label{68}
\end{equation}

In order to simplify (\ref{68}), we denote $\lambda = \frac 1{rv}$ in place of $v$. At this time, $\min_{v\ge 0}$ is replaced by $\max_{\lambda\ge 0}$ and 
\begin{equation}
\frac{r_{1*}^2(\alpha rv+1)^2+r^2v^2\sigma^2+\sigma_*^2}{2(1+rv)rv} = \frac{(\alpha+\lambda)^2r_{1*}^2+\sigma^2+\sigma_*^2\lambda^2}{2(\lambda+1)}. \label{68.1}
\end{equation}

Substituting (\ref{68.1}) in (\ref{67})  we have the final optimization for LASSO:
\begin{equation}
\begin{aligned}
\min_{\alpha\in\R; \atop \sigma\ge0}
\max_{r,\tau,\lambda\ge0; \atop\theta\in \R}
\quad &-\frac{\sigma}{2\tau} -\frac{\alpha^2r_{1*}^2}{2\sigma\tau} 
- \frac{\sigma\tau r^2 \kappa_*}{2} + \frac{r^2\lambda}{2}
+\alpha\theta r_{1*}^2 - \frac{\sigma\tau\theta^2 r_{1*}^2}{2} -\frac{(r_{1*}^2+\sigma_*^2)\lambda}{2} - r_{1*}^2\alpha\\
&+ \frac{(\alpha+\lambda)^2r_{1*}^2+\sigma^2+\sigma_*^2\lambda^2}{2(\lambda+1)} + \mathbb{E}\left[ M_{\lambda_* \tilde{f}}(rZ+\theta b_0;1)\right]\cdot {\sigma\tau \kappa_*}, 
\end{aligned} 
 \label{70}
\end{equation}
which is a smooth function with respect to $\alpha,\sigma,r,\tau,\lambda, \theta $. Let $\phi$ denote the objective function of (\ref{70}).

{\bf Deriving SEs from function $\phi$}: 
The SEs are given by the first order optimality conditions of $\phi$ :
\begin{itemize}

\item[1)]For $\frac{\partial \phi}{\partial \alpha} = 0$: 
\begin{equation}
-\frac{\alpha}{\sigma\tau}+\theta-1+\frac{\alpha+\lambda}{\lambda+1}=0.\label{71.a}
\end{equation}

\item[2)]For $\frac{\partial \phi}{\partial \sigma} = 0$:
\begin{equation}
-\frac{1}{2\tau}-\frac{\tau r^2\kappa_*}{2}+\frac{r_{1*}^2\alpha^2}{2\sigma^2\tau}-\frac{r_{1*}^2\tau\theta^2}{2}+{\tau \kappa_*} E[M_{\lambda_*\tilde{f}}(rZ+\theta b_0;1)]+\frac{\sigma}{\lambda+1}=0. \label{71.b1}
\end{equation}

\item[3)]For $\frac{\partial \phi}{\partial \lambda} = 0$: 
\begin{equation}
r^2-r_{1*}^2-\sigma_*^2+\frac{2[(\alpha+\lambda)r_{1*}^2+\lambda\sigma_*^2]}{\lambda+1}-\frac{(\alpha+\lambda)^2r_{1*}^2+\sigma^2+\lambda^2\sigma_*^2}{(\lambda+1)^2}=0.\label{71.c}
\end{equation}

\item[4)]For $\frac{\partial\phi}{\partial\theta}=0$: 
\begin{equation}
r_{1*}^2\alpha-{\sigma\tau\kappa_*}\mathbb{E}[b_0(Prox_{\lambda_*\tilde{f}}(rZ+\theta b_0))]=0. \label{71.d}    
\end{equation}
where we use the definition $\mathbb{E}[b_0^2] = E_\Pi X^2 = r_{1*}^2$.

\item[5)]For $\frac{\partial\phi}{\partial r}=0$:
\begin{equation}
-{\sigma\tau r\kappa_*}+r\lambda+{\sigma\tau\kappa_*}\mathbb{E}[(rZ+\theta b_0-Prox_{\lambda_*\tilde{f}}(rZ+\theta b_0;1))Z]=0. \label{71.e1}
\end{equation}

Since $\mathbb{E}[Z^2] = 1, \mathbb{E}[Zb_0] = 0$. For any function $\tilde{f}(x)$, the Moreau envelope and proximal operator of $\tilde{f}$ stratifies
\begin{equation}
\frac{\partial}{\partial x} M_{\tilde{f}}(x;t) = \frac {x-Prox_{\tilde f}(x;t)}{t}.
\end{equation}
Then we rewrite the equation (\ref{71.e1}) as
\begin{equation} \label{71.e2}
r\lambda-{\sigma\tau\kappa_*}\mathbb{E}[Prox_{\lambda_*\tilde{f}}(rZ+\theta b_0;1)Z]=0.
\end{equation} 

Using Stein lemma,  
\begin{equation}
\begin{aligned}
\mathbb{E}[Prox_{\lambda_*\tilde{f}}(rZ+\theta b_0;1)Z]&=\mathbb{E}[r\frac{\partial Prox_{\lambda_*\tilde{f}}(rZ+\theta b_0;1)}{\partial x}]\\
&=\mathbb{E}[\frac{\partial Prox_{\lambda_*\tilde{f}}(rZ+\theta b_0;1)}{\partial Z}],
\end{aligned}
\end{equation}

(\ref{71.e2}) can be rewritten as
\begin{equation}
   \lambda={\sigma\tau\kappa_*}\mathbb{E}[\frac{\partial Prox_{\lambda_*\tilde{f}}(rZ+\theta b_0;1))}{\partial x}].\label{71.e} 
\end{equation}

\item[6)]For $\frac{\partial\phi}{\partial \tau}=0$: 
\begin{equation}
\frac{\sigma}{2\tau^2}-\frac{\sigma r^2\kappa_*}{2}+\frac{r_{1*}^2\alpha^2}{2\sigma\tau^2}-\frac{r_{1*}^2\sigma\theta^2}{2}+{\sigma\kappa_*} \mathbb{E}[M_{\lambda_*\tilde{f}}(rZ+\theta b_0;1)]=0. \label{71.f1}
\end{equation}

For any function $\tilde{f}(x)$, the Moreau envelope and proximal operator of $\tilde{f}$ stratifies
\begin{equation}
M_{\lambda_*\tilde{f}}(x;b)=\lambda_*M_{\tilde{f}}(x;\lambda_*b)=\frac{x^2}{2b}-\frac{[Prox_{\tilde{f}}(x;\lambda_*b)]^2}{2b}, \quad \forall \lambda_*,b>0, x\in \R.
\end{equation}

Using this property, we can rewrite (\ref{71.f1}) as:
\begin{equation}
\frac{\sigma}{2\tau^2}-\frac{\sigma r^2\kappa_*}{2}+\frac{r_{1*}^2\alpha^2}{2\sigma\tau^2}-\frac{r_{1*}^2\sigma\theta^2}{2}+{\sigma\kappa_*} \mathbb{E}\left\{\frac{(rZ+\theta b_0)^2}{2}-\frac{[Prox_{\lambda_*\tilde{f}}(rZ+\theta b_0;1)]^2}{2}\right\}=0.\\
\end{equation}

i.e., 
\begin{equation}
    \frac{\sigma}{2\tau^2}+\frac{r_{1*}^2\alpha^2}{2\sigma\tau^2}-\frac{\sigma\kappa_*}{2} \mathbb E[(Prox_{\lambda_*\tilde{f}}(rZ+\theta b_0;1))^2]=0.\label{71.f}
\end{equation}
\end{itemize}
Similarly, the equation (\ref{71.b1}) derived by $\frac{\partial \phi}{\partial \sigma} = 0$ can be rewritten as
\begin{equation}
-\frac{1}{2\tau}-\frac{\tau r^2\kappa_*}{2}+\frac{r_{1*}^2\alpha^2}{2\sigma^2\tau}-\frac{r_{1*}^2\tau\theta^2}{2}+{\tau\kappa_*} E[\frac{(rZ+\theta b_0)^2}{2}-\frac{[Prox_{\lambda_*\tilde{f}}(rZ+\theta b_0;1)]^2}{2}]+\frac{\sigma}{\lambda+1}=0
\end{equation}
i.e.,
\begin{equation}
    -\frac{1}{2\tau}+\frac{r_{1*}^2\alpha^2}{2\sigma^2\tau}-\frac{\tau\kappa_*}{2} E[(Prox_{\lambda_*\tilde{f}}(rZ+\theta b_0;1))^2]+\frac{\sigma}{\lambda+1}=0.\label{71.b}
\end{equation}

Hence we get the SEs by summarizing equations (\ref{71.a}),  (\ref{71.b}),  (\ref{71.c}),  (\ref{71.d}),  (\ref{71.e}), (\ref{71.f}):
\begin{equation} \label{72}
\begin{aligned}
0&=-\frac{\alpha}{\sigma\tau}+\theta-1+\frac{\alpha+\lambda}{\lambda+1},\\
0&=-\frac{1}{2\tau}+\frac{r_{1*}^2\alpha^2}{2\sigma^2\tau}-\frac{\tau\kappa_*}{2} E[(Prox_{\lambda_*\tilde{f}}(rZ+\theta b_0;1))^2]+\frac{\sigma}{\lambda+1},\\
0&=r^2-r_{1*}^2-\sigma_*^2+\frac{2[(\alpha+\lambda)r_{1*}^2+\lambda\sigma_*^2]}{\lambda+1}-\frac{(\alpha+\lambda)^2r_{1*}^2+\sigma^2+\lambda^2\sigma_*^2}{(\lambda+1)^2},\\
0&=r_{1*}^2\alpha-{\sigma\tau\kappa_*}E[b_0(Prox_{\lambda_*\tilde{f}}(rZ+\theta b_0;1))],\\
\lambda&= {\sigma\tau\kappa_*}E[\frac{\partial Prox_{\lambda_*\tilde{f}}(rZ+\theta b_0;1))}{\partial x}],\\
0&=\frac{\sigma}{2\tau^2}+\frac{r_{1*}^2\alpha^2}{2\sigma\tau^2}-\frac{\sigma\kappa_*}{2} E[(Prox_{\lambda_*\tilde{f}}(rZ+\theta b_0;1))^2].
\end{aligned}
\end{equation}
regarding $(\alpha,\sigma, \lambda, \theta, r, \tau)$.

Since $r_*^2 = \mathbb{E}_{b_0\sim \Pi}b_0^2 = \frac{r_{1*}^2}{\kappa_*}$, the SEs (\ref{72}) can be rewritten as
\begin{subequations} \label{FinalSEofCGMT}
\begin{align}
0&=-\frac{\alpha}{\sigma\tau}+\theta-1+\frac{\alpha+\lambda}{\lambda+1}\label{71a},\\
0&=-\frac{1}{2\tau}+\frac{r_*^2 \kappa_*\alpha^2}{2\sigma^2\tau}-\frac{\tau\kappa_*}{2} E[(Prox_{\lambda_*\tilde{f}}(rZ+\theta b_0;1))^2]+\frac{\sigma}{\lambda+1}\label{71b},\\
0&=r^2-r_*^2 \kappa_*-\sigma_*^2+\frac{2[(\alpha+\lambda)r_*^2 \kappa_*+\lambda\sigma_*^2]}{\lambda+1}-\frac{(\alpha+\lambda)^2r_*^2 \kappa_*+\sigma^2+\lambda^2\sigma_*^2}{(\lambda+1)^2}\label{71c},\\
0&=r_*^2 \kappa_*\alpha-{\sigma\tau\kappa_*}E[b_0(Prox_{\lambda_*\tilde{f}}(rZ+\theta b_0;1))]\label{71d},\\
\lambda&= {\sigma\tau\kappa_*}E[\frac{\partial Prox_{\lambda_*\tilde{f}}(rZ+\theta b_0;1))}{\partial x}]\label{71e},\\
0&=\frac{\sigma}{2\tau^2}+\frac{r_*^2 \kappa_*\alpha^2}{2\sigma\tau^2}-\frac{\sigma\kappa_*}{2} E[(Prox_{\lambda_*\tilde{f}}(rZ+\theta b_0;1))^2].\label{71f}
\end{align}    
\end{subequations}
regarding $(\alpha,\sigma, \lambda, \theta, r, \tau)$. This completes the proof.

\section{Proof of Theorem 2} \label{appendixC.1:Equivalence of SEs}

We first rewrite $r$, $\tau$, $Z$ and $b_0$ in (\ref{FinalSEofCGMT}) to $\gamma_2$, $\tau_2$ and $Z_2$ $\beta_2$ respectively, the equation (\ref{71e}) becomes
\begin{equation}
    \lambda =\sigma\tau_2\kappa_* \mathbb E[\frac{\partial Prox_{\tilde{f}}(\gamma_2Z_2+\theta \beta_2;\lambda_*)}{\partial x}] = \sigma\tau_2\kappa_* \mathbb E[\eta'(\gamma_2Z_2+\theta \beta_2;\lambda_*)]
\end{equation}
for $\tilde{f}(x) = |x|$.

Then we simplify the SEs (\ref{FinalSEofCGMT}). Consider equations (\ref{71b}) and (\ref{71f}) and we have 
\begin{equation} \label{SEequiv1}
    \sigma\tau_2 = \lambda + 1,
\end{equation}
substituting (\ref{SEequiv1}) in (\ref{71a}) we have 
\begin{equation}\label{73} 
    \theta=\frac{1}{\lambda+1}=\frac{1}{\sigma\tau_2}.
\end{equation}
For the second equation of AMP:
\begin{equation}
    \gamma_1 =\kappa_*(\gamma_1+\lambda_*)\mathbb E[\eta'(\beta_1+\tau_1 Z_1;\lambda_*+\gamma_1)].
\end{equation}
it is obviously equivalent to the equation (\ref{71e}) if we have parameter transformations $\tau_1 = \frac{\gamma_2}{\theta} $ and $\gamma_1 = \lambda_*(\frac 1\theta - 1)$. A property of $\eta(\cdot;\cdot)$ is used for the equivalence:
\begin{equation}
    \eta'(cx;ct) = \eta'(x;t), \quad \forall c>0.
\end{equation}

Using the parameter transformations mentioned above and denote $W = \eta(\beta_1+\tau_1Z_1;\lambda_*+\gamma_1)$, the first equation of AMP in Proposition 3.1:
\begin{equation}\label{the first equation of AMP in Proposition 3.1}
\tau_1^2=\sigma_*^2+\kappa_* \mathbb E[\eta(\beta_1+\tau_1 Z_1;\lambda_*+\gamma_1)-\beta_1]^2
\end{equation}
can be rewritten as 
\begin{equation}
    \frac{\gamma_2^2}{\theta^2}=\sigma_*^2+\kappa_* \mathbb{E}[W-\beta_2]^2,
\end{equation}
i.e., 
\begin{equation} \label{79}
    \frac{\gamma_2^2}{\theta^2}=\sigma_*^2+\kappa_*(\mathbb{E}[W^2]+\mathbb{E}[\beta_2^2]-2\mathbb{E}[W\beta_2]).
\end{equation}
Substituting (\ref{71d}), (\ref{71f}) and $\mathbb{E}(\beta_2^2)=r_*^2$, we have
\begin{equation}
  \begin{aligned}
  \mathbb{E}(W^2) &= \mathbb{E}(\eta^2(\beta_1+\tau_1Z_1;\lambda_*+\gamma_1))\\
  &= \mathbb{E}(\eta^2(\beta_1+\frac{\gamma_2}{\theta}Z_1;\frac{\lambda_*}{\theta}))\\
  &= \frac{1}{\theta^2}\mathbb{E}(\eta^2(\theta\beta_2+\gamma_2 Z_2; \lambda_*)) \quad(\text{because }\eta(cx;ct) = c\eta(x;t))\\
  &= \frac{1}{\theta^2}\left[\frac{1}{\tau_2^2\kappa_*}+ \frac{r_*^2\alpha^2}{\sigma^2\tau_2^2} \right] \quad(\text{using } (\ref{71f})), \\
  -2\mathbb{E}(W\beta_2) &= -2\mathbb{E}(\eta(\beta_1+\tau_1Z_1;\lambda_*+\gamma_1)\cdot \beta_1)\\
  &= -\frac{2}{\theta} \mathbb{E}(\eta(\theta\beta_2+\gamma_2 Z_2; \lambda_*)\cdot \beta_2) \\
  &= -\frac{2}{\theta}\frac{r_*^2\alpha}{\sigma\tau_2} \quad(\text{using } (\ref{71d})). \\
  \end{aligned}  
\end{equation}
Then (\ref{79}) can be rewritten as
\begin{equation}\label{80}
\begin{aligned}
    \frac{r^2}{\theta^2}&=\sigma_*^2+\kappa_*\left(\frac{1}{\theta^2\tau_2^2\kappa_*}+ \frac{r_*^2\alpha^2}{\theta^2\sigma^2\tau_2^2} - \frac{2r_*^2\alpha}{\theta\sigma\tau_2 + r_*^2} \right).
\end{aligned}
\end{equation}
Using (\ref{73}) in (\ref{80}) we have
\begin{equation} \label{81}
    r^2(\sigma\tau_2)^2 = \sigma_*^2 + \sigma^2 +  \kappa_*(\alpha-1)^2r_*^2.
\end{equation}
Besides, for CGMT, the equation (\ref{71c}) can be written as 
\begin{equation}\label{82}
    (\lambda+1)^2r^2-(\alpha-1)^2r_*^2\kappa_*-\sigma^2-\sigma_*^2=0.
\end{equation} 
Using (\ref{73}) in (\ref{82}) we have
\begin{equation}\label{83}
    (\sigma\tau_2)^2r^2=\sigma_*^2+(\alpha-1)^2r_*^2\kappa_*+\sigma^2.
\end{equation}
 The equation (\ref{81}) and (\ref{83}) are equivalent. Hence the equations (\ref{71d}), (\ref{71f}) and (\ref{71c}) of CGMT can be shown to be a decomposition of the equation (\ref{the first equation of AMP in Proposition 3.1}) after some parameter transformations. In conclusion we prove the equivalence between SEs from CGMT and AMP in LASSO framework.

\section{Proof of Theorem 3}\label{appendix: equivalence of SEs of logistic regression from LOO and CGMT}
By doing the following parameter transformations:
$$
\alpha_2=\sqrt{\kappa_*}\alpha_1, \mu=r_*\sigma, \lambda_2=\lambda_1,
$$
SEs of CGMT become:
\begin{equation}
  \begin{aligned}
0&=\mathbb E[Vl'(Prox_{\lambda_1l}(\sqrt{\kappa_*}\alpha_1 Z+r_*\sigma V))],\\
\kappa_*^2(\alpha_1)^2&=(\lambda_1)^2\mathbb E[(l'(Prox_{\lambda_1 l}(\sqrt{\kappa_*}\alpha_1 Z+r_*\sigma V)))^2],\\
\kappa_*&=\lambda_1\mathbb E[\frac{l''(Prox_{\lambda_1l}(\sqrt{\kappa_*}\alpha_1 Z+r_*\sigma V))}{1+\lambda_1l''(Prox_{\lambda_1l}(\sqrt{\kappa_*}\alpha_1 Z+r_*\sigma V))}].
\end{aligned} 
\end{equation}
What we want to prove is that:
$$
\begin{aligned}
\mathbb E[(l'(prox_{\lambda_1 l}(\sqrt{\kappa_*}\alpha_1 Z+r_*\sigma V)))^2]&=\mathbb E[2\rho'(Q_1)(\rho'(prox_{\lambda_1\rho}(Q_2)))^2],\\
1-\lambda_1\mathbb E[\frac{l''(prox_{\lambda_1l}(\sqrt{\kappa_*}\alpha_1 Z+r_*\sigma V))}{1+\lambda_1l''(prox_{\lambda_1l}(\sqrt{\kappa_*}\alpha_1 Z+r_*\sigma V))}]&=\mathbb E[\frac{2\rho'(Q_1)}{1+\lambda_1\rho''(prox_{\lambda_1\rho}(Q_2))}],\\
E[V l'(prox_{\lambda_1l}(\sqrt{\kappa_*}\alpha_1 Z+r_*\sigma V))]&=c\mb E[\rho'(Q_1)Q_1\rho'(prox_{\lambda_1\rho}(Q_2))].
\end{aligned}
$$
where $c$ is a constant.

First, we verify the following identity:
\begin{equation}\label{equation, aim1}
 \mathbb E[(l'(prox_{\lambda_1 l}(\sqrt{\kappa_*}\alpha_1 Z+r_*\sigma V)))^2]=\mathbb E[2\rho'(Q_1)(\rho'(prox_{\lambda_1\rho}(Q_2)))^2].   
\end{equation}
Note that
\begin{equation}
\begin{aligned}
l(t)&=\rho(-t),\\l'(t)&=-\rho'(-t),\\
l''(t)&=\rho(t),\\Prox_{\lambda_1l}(z)&=-Prox_{\lambda_1\rho}(-z).
\end{aligned}
\end{equation}
and the probability density function (pdf) of $V=GY$ is $$P_V(v)=\frac{1}{\sqrt{2\pi}}e^{-\frac{v^2}{2}}\frac{2}{1+e^{-r_* v}},$$
we have:
\begin{equation}
\begin{aligned}
&\mathbb E[(l'(Prox_{\lambda_1 l}(\sqrt{\kappa_*}\alpha_1 Z+r_*\sigma V)))^2]\\
=&\iint(l'(Prox_{\lambda_1 l}(\sqrt{\kappa_*}\alpha_1 h+r_*\sigma v)))^2P_Z(h)P_V(v)dhdv\\
=&\iint2(\rho'(Prox_{\lambda_1\rho}(-\sqrt{\kappa_*}\alpha_1 h-r_*\sigma v)))^2\frac{1}{2\pi}e^{-\frac{h^2+v^2}{2}}\rho'(r_* v)dhdv
\end{aligned}
\end{equation}
where $P_Z(h):=\frac{1}{\sqrt{2\pi}}e^{-\frac{h^2}{2}}$ is the pdf of $Z$.
Meanwhile,
\begin{equation}\label{equation, Q1,Q2}
\begin{aligned}
\mathbb E[2\rho'(Q_1)(\rho'(prox_{\lambda_1\rho}(Q_2)))^2]&=\iint2\rho'(q_1)(\rho'(prox_{\lambda_1\rho}(q_2)))^2P_{Q_1,Q_2}(q_1,q_2)dq_1dq_2.\\
\end{aligned}
\end{equation}

Now we introduce the following parameter transformations: $q_1=r_* v,q_2=-\sqrt{\kappa_*}\alpha_1 h-r_*\sigma v$. Then (\ref{equation, Q1,Q2}) becomes
$$
\iint\sqrt{\kappa_*}\alpha_1 r_**2\rho'(r_* v)(\rho'(prox_{\lambda_1\rho}(-\sqrt{\kappa_*}\alpha_1 h-r_*\sigma v)))^2P_{Q_1,Q_2}(r_* v,-\sqrt{\kappa_*}\alpha_1 h-r_*\sigma v)dhdv.
$$

In order to verify (\ref{equation, aim1}), we only need to prove that:
$$
\frac{1}{2\pi}e^{-\frac{h^2+v^2}{2}}=\sqrt{\kappa_*}\alpha_1 r_* P_{Q_1,Q_2}(r_* v,-\sqrt{\kappa_*}\alpha_1 h-r_*\sigma v).
$$
Construct $Q_1',Q_2'$ as follows: assume $Z',V'\overset{i.i.d.}{\sim}\mathcal N(0,1)$ and 
$$
\left(
\begin{matrix}
   Q_1'  \\
   Q_2'
  \end{matrix} 
  \right)
=
\left(
\begin{matrix}
   0 & r_*\\
   -\sqrt{\kappa_*}\alpha_1 & -r_*\sigma
  \end{matrix} 
  \right)
  \left(
\begin{matrix}
   Z'  \\
   V'
  \end{matrix} 
  \right).
$$
We can easily verify that:
$$
\begin{aligned}
\mathbb E[(Q_1',Q_2')^T]&=(0,0)^T ,\\
Cov[(Q_1',Q_2')^T]
&=\left(
\begin{matrix}
   r_*^2 & -r_*\sigma r_*  \\
   -r_*\sigma r_* & \sqrt{\kappa_*}\alpha_1^2+r_*\sigma^2 
  \end{matrix}
\right)
\end{aligned}
$$
which means $(Q_1',Q_2')$ has identical distribution of $(Q_1,Q_2)$.

On the other hand, since
$$
\begin{aligned}
P_{Q_1',Q_2'}(q_1',q_2')d q_1'q_2'&=P_{Z',V'}(h',v')d h' d v', \\
\frac{d q_1'q_2'}{dh'dv'}&=r_*\sqrt{\kappa_*}\alpha_1,\\ q_1'&=r_* v',\\q_2'&=-\sqrt{\kappa_*}\alpha_1 h'-r_*\sigma v',
\end{aligned}
$$
we have:
$$
\frac{1}{2\pi}e^{-\frac{h'^2+v'^2}{2}}=\sqrt{\kappa_*}\alpha_1 r_* P_{Q_1',Q_2'}(r_* v',-\sqrt{\kappa_*}\alpha_1 h'-r_*\sigma v'),
$$
which completes our proof of (\ref{equation, aim1}).

Secondly, we prove:
\begin{equation}\label{equation, aim2}
1-\lambda_1\mathbb E[\frac{l''(prox_{\lambda_1l}(\sqrt{\kappa_*}\alpha_1 Z+r_*\sigma V))}{1+\lambda_1l''(prox_{\lambda_1l}(\sqrt{\kappa_*}\alpha_1 Z+r_*\sigma V))}]=\mathbb E[\frac{2\rho'(Q_1)}{1+\lambda_1\rho''(prox_{\lambda_1\rho}(Q_2))}].\\
\end{equation}

Left hand side (LHS) of (\ref{equation, aim2}) is
\begin{equation}
\begin{aligned}
&\mb E[\frac{1}{1+\lambda_1l''(prox_{\lambda_1l}(\sqrt{\kappa_*}\alpha_1 Z+r_*\sigma V))}]\\
=&\iint\frac{1}{1+\lambda_1l''(prox_{\lambda_1l}(\sqrt{\kappa_*}\alpha_1 h+r_*\sigma v))}P_Z(h)P_V(v)dhdv\\
=&\iint\frac{1}{1+\lambda_1l''(prox_{\lambda_1l}(\sqrt{\kappa_*}\alpha_1 h+r_*\sigma v))}\frac{1}{2\pi}e^{-\frac{h^2+v^2}{2}}2\rho'(r_* v)dhdv\\
=&\iint\frac{1}{1+\lambda_1\rho''(prox_{\lambda_1\rho}(-\sqrt{\kappa_*}\alpha_1 h-r_*\sigma v))}\frac{1}{2\pi}e^{-\frac{h^2+v^2}{2}}2\rho'(r_* v)dhdv.
\end{aligned}    
\end{equation}

Through parameter transformations $q_1=r_* v,q_2=-\sqrt{\kappa_*}\alpha_1 h-r_*\sigma v$, RHS of (\ref{equation, aim2}) becomes 
\begin{equation}
\begin{aligned}
&\iint\frac{2\rho'(q_1)}{1+\lambda_2\rho''(prox_{\lambda_1\rho}(q_2))}P_{Q_1,Q_2}(q_1,q_2)dq_1dq_2\\=
&\iint\frac{2\rho'(r_* v)}{1+\lambda_1\rho''(prox_{\lambda_1\rho}(-\sqrt{\kappa_*}\alpha_1 h-r_*\sigma v))}P_{Q_1,Q_2}(r_* v,-\sqrt{\kappa_*}\alpha_1 h-r_*\sigma v)r_*\sqrt{\kappa_*}\alpha_1 dhdv.
\end{aligned}
\end{equation}

Combining with
$$
\frac{1}{2\pi}e^{-\frac{h^2+v^2}{2}}=P_{Q_1,Q_2}(r_* v,-\sqrt{\kappa_*}\alpha_1 h-r_*\sigma v)r_*\sqrt{\kappa_*}\alpha_1
$$
completes the proof of (\ref{equation, aim2}).

The proof of
$$
\mb E[V l'(prox_{\lambda_1l}(\sqrt{\kappa_*}\alpha_1 Z+r_*\sigma V))]=c\mb E[\rho'(Q_1)Q_1\rho'(prox_{\lambda_1\rho}(Q_2))]
$$
can be derived similarly. So we omit the proof here.

\section{relaxation phenomenon}
\subsection{relaxation phenomenon of M-estimator}\label{Appendix, relaxation phenomenon of M-estimator}

Notice that (\ref{relax-M-latter}) is equivalent to 
\begin{equation}\label{relax-M-end1}
  \min_{\tau_3,\boldsymbol v}\frac{1}{n}\sum_{i=1}^n\rho(v_i)\quad s.t.\frac{1}{\sqrt{n}}||\tau_3\boldsymbol h+\boldsymbol v-\bm\epsilon||_2\leq\frac{\tau_3}{\sqrt{n}}||\boldsymbol g||_2.
\end{equation}
On the other hand, by rotation invariance of Gaussian distribution, (\ref{relax-M-former}) becomes:
\begin{equation}\label{relax-M-end2}
    \min_{||\boldsymbol w||_2,\boldsymbol v}\frac{1}{n}\sum_{i=1}^n\rho(v_i)\quad s.t.\ ||\boldsymbol v-\boldsymbol  \epsilon+||\boldsymbol w||_2Z'||=0.
\end{equation}
where $Z'\sim\mathcal N(0,1),\boldsymbol\epsilon:=(\epsilon_1,...,\epsilon_n)^T$.
Comparing (\ref{relax-M-end1}) with (\ref{relax-M-end2}) can verify the relaxation phenomenon, i.e. the only difference between AO and PO is that the feasible region of PO is a subset of the feasible region of AO.

\subsection{relaxation phenomenon of support vector machine and Logistic regression}\label{Appendix, relaxation phenomenon of support vector machine and Logistic regression}
The relaxation phenomenon of support vector machine and logistic regression can be similarly shown as we have done in Appendix \ref{Appendix, relaxation phenomenon of M-estimator}, so we omit the proof here.

\bibliography{reference}
\bibliographystyle{chicago}